\documentclass[12pt, centertags]{amsart}
\usepackage{amsmath, amstext, amsthm, a4, amssymb, amscd}
\usepackage{amssymb, graphics, epsfig, color}
\usepackage[mathscr]{eucal}
\usepackage[all]{xy}
\usepackage{mathrsfs}
\usepackage{graphicx}
\makeindex

\numberwithin{equation}{section}

\newcommand{\etalchar}[1]{$^{#1}$}

 \newcommand{\field}[1]{\mathbb{#1}}
 \newcommand{\Z}{\field{Z}} \newcommand{\R}{\field{R}} \newcommand{\C}{\field{C}} \newcommand{\N}{\field{N}} 
 
 \def\cC{\mathscr{C}}

    \def\Re{{\rm Re}}

 \DeclareMathOperator{\Ker}{Ker} \DeclareMathOperator{\Coker}{Coker} \DeclareMathOperator{\Dom}{Dom}
 \DeclareMathOperator{\ran}{Range}  \DeclareMathOperator{\Id}{Id} 
 \DeclareMathOperator{\tr}{Tr} \DeclareMathOperator{\str}{Tr_{s}}

 \DeclareMathOperator{\End}{End}

\DeclareMathOperator{\sign}{sign}\DeclareMathOperator{\ind}{Ind}
\DeclareMathOperator{\erfc}{erfc}

 \newcommand{\norm}[1]{\lVert#1\rVert} \newcommand{\abs}[1]{\lvert#1\rvert} 
 
 \newcommand{\iprod}[1]{\langle #1\rangle }
 \newcommand{\pd}[2]{\frac{\partial #2}{\partial #1}}

\newcommand{\cY}[1]{Y_{[#1]}}

\newcommand{\bc}{\mathbf{c}}
\newcommand{\nS}{\nabla^S}

\newcommand{\Domm}{(D_{M_1}+mf)_{V_-}}
\newcommand{\Dtmp}{(D_{M_2}+mf)_{V_+}}

 \newtheorem{thm}{Theorem}[section] \newtheorem{lemma}[thm]{Lemma}  \newtheorem{prop}[thm]{Proposition} 
 
 \newtheorem{cor}[thm]{Corollary} \theoremstyle{definition} \newtheorem{rem}[thm]{Remark} \theoremstyle{definition}
 \newtheorem{defn}[thm]{Definition}
 \newtheorem{Assumption}[thm]{Assumption} %
 \theoremstyle{remark} %
 \newcommand{\be}{\begin{eqnarray}}
   \newcommand{\var}{\varepsilon}
 \newcommand{\comment}[1]{}
\newcommand{\onehalf}{\frac{1}{2}}
\newcommand{\Cinf}[1]{C^{\infty}\left( #1\right) }
\newcommand{\Cinfc}[1]{C^{\infty}_0\left( #1\right) }

\newcommand{\ra}{\rightarrow}

\newcommand{\rac}[1]{\overset{#1}\ra}

\newcommand{\n}{\nabla}
\newcommand{\set}[1]{\{#1\}}
\newcommand{\GS}{\Gamma_S}

\newcommand{\tM}{\widetilde{M}}

\def\cE{\mathscr{E}}

\newcommand{\pr}[1]{{\rm pr}_{#1}}

\newcommand{\tPi}{\widetilde{\Pi}}
\newcommand{\CLM}{{\rm Cl}(M)}

\newcommand{\G}{\Gamma}

\newcommand{\gTY}{g^{TY}}

\newcommand{\Ymoo}{Y_{[-1,1]}}

\newcommand{\Sp}{S_{+}}\newcommand{\Sm}{S_{-}}

\newcommand{\LtME}{L^2(M,S)}\newcommand{\GME}{\Cinf{M,S}}
\newcommand{\lj}{\lambda_j}

\newcommand{\gTM}{g^{TM}}\newcommand{\Yzo}{Y_{[0,1]}}\newcommand{\Ymo}{Y_{[-1,0]}}
\newcommand{\Mp}{M_{+}}\newcommand{\Mm}{M_{-}}

\newcommand{\bM}{\overline{M}}
\newcommand{\DDWm}{D_{\tiny DW,m}}\newcommand{\DPVm}{D_{\tiny PV,m}}

\newcommand{\Tmo}{\frac{1}{T}}

\begin{document}

 \title{Atiyah-Patodi-Singer index and Domain-wall eta invariants}
 \date{\today}
\author{Jialin Zhu}
 \address{Mathematical Science Research Center, Chongqing University of Technology, No. 69 Hongguang Road, Chongqing 400054, China}
 \email{leozjl@mail.ustc.edu.cn}

\begin{abstract}

In this paper we establish a formula, expressing the generalized Atiyah-Patodi-Singer index in terms of  eta invariants of domain-wall massive Dirac operators, without assuming that the Dirac operator on the boundary is invertible. Compared with the original Atiyah-Patodi-Singer index theorem, this formula has the advantage that no global spectral projection boundary conditions appear. Our main tool is an asymptotic gluing formula for eta invariants proved by using a splitting principle developed by Douglas and Wojciechowski in adiabatic limit. The eta invariant splits into a contribution from the interior, one from the boundary, and an error term vanishing in the adiabatic limit process.

\end{abstract}

\maketitle

\tableofcontents
\setcounter{section}{-1}

\section{Introduction}\label{s0}

In this paper we establish a formula expressing the generalized Atiyah-Patodi-Singer index in terms of the eta invariants of domain-wall massive Dirac operators without assuming that the Dirac operator on the boundary is invertible. Our results generalize the work of H. Fukaya, M. Furuta, S. Matsuo, T. Onogi, S. Yamaguchi and  M. Yamashita \cite{FF20}. They are motivated by the study of lattice gauge theory to investigate such a problem.

The Atiyah-Patodi-Singer index was introduced by Atiyah, Patodi and Singer \cite{APS} in their generalization of the Atiyah-Singer index theorem to manifolds with boundary. In \cite{APS}, they also introduced a global boundary condition by using the spectral projection operator of the Dirac operator on the boundary. An important spectral invariant, called eta invariant, appears in their index formula as a correction term from the boundary.  The Atiyah-Patodi-Singer index and the eta invariant have important applications in condensed matter physics, especially in studying the symmetry-protected topological phases of matter (cf. \cite{Wi16}, \cite{Yon16}, \cite{OY21}). In symmetry-protected topological phases of matter, physicists use massive fermion Dirac operators and impose local boundary conditions, while massless Dirac operators and global boundary conditions are used in the Atiyah-Patodi-Singer index theorem. The global spectral projection boundary condition is physicist-unfriendly, since any boundary conditions on fields must be local in relativistic physics otherwise the information will propagate faster than the speed of light along the boundary  (cf. \cite{Fuk21}). In \cite{APS}, Atiyah, Patodi and Singer initially tried to use local boundary conditions, such as the Dirichlet or Neumann boundary conditions, but there are topological obstructions if one wants to make the index of signature operator and the signature of the manifold coincide.

In \cite{FOY17}, H. Fukaya, T. Onogi and S. Yamaguchi first derived a formula, relating the Atiyah-Patodi-Singer index to the eta invariant of domain-wall fermion Dirac operator, for a four dimensional flat manifold. In \cite{FF20}, H. Fukaya, M. Furuta, S. Matsuo, T. Onogi, S. Yamaguchi and  M. Yamashita proved this formula for any even-dimensional Riemannnian manifolds by using an embedding trick and a Witten localisation argument (cf. \cite{Wit82}), on the assumption that the boundary Dirac operator is invertible. In this paper we prove that the formula still holds even if the Dirac operator on the boundary is not invertible, based on the adiabatic splitting principle developed by Douglas and Wojciechowski (cf. \cite{DouWoj91}). It is natural to expect that there is a new term coming from the boundary in our formula. This new term is related to the choice of Lagrangian subspaces in the kernel of the Dirac operator on the boundary. The Lagrangian subspaces are used in the definition of the generalized Atiyah-Patodi-Singer boundary conditions. This term will disappear if the Dirac operator on the boundary is invertible.

In \cite{DouWoj91}, Douglas and Wojciechowski first studied the adiabatic limit of eta invariants with the generalized Atiyah-Patodi-Singer boundary condition on manifolds with boundary. They developed a general method for calculating the adiabatic limit of spectral invariants. In \cite{DouWoj91}, the product structures of metrics are assumed on a cylinder neighborhood of the boundary, and the adiabatic limit is a process of making the length of this cylinder go to infinity. They found that the eta invariant associated with compatible Dirac operators is split into a contribution from the interior, one from the cylinder and an error term vanishing with
the increase of the length of the cylinder, on the assumption that the boundary Dirac operator is invertible. In \cite{Muller94}, M\"uller studied the behavior of eta invariants with generalized Atiyah-Patodi-Singer boundary conditions in the adiabatic limit when the Dirac operator on the boundary is not invertible. When the boundary Dirac operator is invertible, the nonzero eigenvalues of the interior Dirac operator stay bounded away from zero. But if the boundary Dirac operator is not invertible, the situation becomes much more difficult since the nonzero spectrum of the interior Dirac operator will cluster at zero as the length of the cylinder part near the boundary goes to infinity. In \cite{Muller94}, scattering matrices are used to analyze the limiting behavior of these small eigenvalues. In \cite{LW96}, Lesch and Wojciechowski studied how the eta invariants depend on the choice of different generalized Atiyah-Patodi-Singer boundary conditions.

In \cite{Wojcie}, Wojciechowski used their method to prove the gluing formula of eta invariants for compatible Dirac operators, on the assumption that the boundary Dirac operator is invertible.
Later in \cite{Wo95}, he established the gluing formula of eta invariants for compatible Dirac operators in the non-invertible case. In \cite{Wojcie}, \cite{Wo95}, the gluing formula holds modulo some integers since the eta invariant is sensitive to the variation of the metrics and can have integer jumps. In \cite{bunke95}, Bunke proved the gluing formula of eta invariants with generalized Atiyah-Patodi-Singer boundary conditions for Dirac-type operators and without the invertibility of the boundary Dirac operator. The gluing formula in \cite{bunke95} holds in real number field,  and some additional non-explicit integer terms appear in the formula. Bunke proved that these integer terms vanish in the adiabatic limit under certain regularity conditions. In this paper we will prove an asymptotic gluing formula for eta invariants, without modulo some integers, by using the adiabatic limit method. The asymptotic gluing formula allows us to reduce our problem to a small neighborhood of the boundary.

The adiabatic limit method has proved to be a very effective tool in the study of spectral invariants, like eta invariant, analytic torsion \cite{RaySing71}. Following \cite{Muller94} and \cite{PaWo06}, Puchol, Zhang and the author gave a pure analytic proof of the gluing formula of Ray-Singer analytic torsion in \cite{PZZ21} by using the adiabatic limit method.

Now we explain our results in more detail.

Let $M$ be a compact oriented manifold of even dimension $m$ with boundary $Y$, where $Y$ is a compact smooth manifold of dimension $m-1$. Let $S=S_+\oplus S_-$ be a $\Z_2-$graded vector bundle over $M$. We denote the $\Z_2-$grading operator by $\Gamma_S$, i.e., $\Gamma_S |_{S_{\pm}}=\pm 1$.

 Let $\gTM$ be a Riemannian metric on $M$. Let $h^S$ be a Hermitian metric over $S$, such that $S_+$ and
$S_-$ are orthogonal. We assume that both $\gTM$ and $h^S$ have product structures near the boundary $Y$.
Let $Y\times [0,1]$ be a product neighborhood of $Y\cong Y\times\set{0}\subset M$. Let $\gTY$ be the Riemannian metric on $Y$ induced by $\gTM$. We assume that the Riemannian metric $\gTM$ is of product form on $Y\times [0,1]$, i.e.,
\begin{align}\begin{aligned}\label{e.17}
\gTM(y,u)=\gTY(y)+du^2,\quad (y,u)\in Y\times [0,1].
\end{aligned}\end{align}
We assume that
\begin{align}\begin{aligned}\label{e.18}
S|_{Y\times [0,1]}=\psi^*(S|_{Y})\quad \text{and} \quad h^{S}|_{Y\times [0,1]}=\psi^*(h^{S}|_{Y}),
\end{aligned}\end{align}
where $\psi:Y\times [0,1]\ra Y$ is the natural projection.

Let $D_M:\Cinf{M, S}\ra \Cinf{M,S}$ be an elliptic, formally self-adjoint, first order differential operator. We assume that $D_M$ is of odd parity, i.e., $D_M\GS=-\GS D_M$. On the product collar, we assume that
\begin{align}\begin{aligned}\label{e.292}
D_M=\gamma\Big(\pd{u}{\,}+D_Y\Big),
\end{aligned}\end{align}
where $\gamma: S|_Y\ra S|_Y$ is a unitary bundle isomorphism and $D_Y:\Cinf{Y,S|_Y}\ra \Cinf{Y,S|_Y}$ is a self-adjoint first order elliptic operator on $Y$.
Assume that
\begin{align}\begin{aligned}\label{e.293}
&\gamma^2=-\Id,\quad \gamma\Gamma_S=-\Gamma_S\gamma,\quad\gamma^*=-\gamma,\\
&\quad \gamma D_Y=-D_Y\gamma,\quad \GS D_Y=D_Y\GS.
\end{aligned}\end{align}
If the Hermitian vector bundle $S$ has a Clifford-module structure, there is a nature chirality operator $\GS$ over $S$ (see \eqref{e.188}) and the associated Dirac operator naturally satisfies the above assumptions \eqref{e.292}, \eqref{e.293}.

Let $P_{>}$ (resp. $P_{<}$) be the projection onto the subspace of $L^2(Y,S|_Y)$ spanned by the eigenvectors corresponding to positive (resp. negative) eigenvalues of $D_Y$. There is a symplectic structure $\Phi$ on the kernel of $D_Y$ induced by the isomorphism $\gamma$, that is for $x,y\in \Ker D_Y\subset \Cinf{Y,S|_Y}$,
\begin{align}\begin{aligned}\label{e.294}
\Phi(x,y):=\iprod{\gamma x,y}_Y,
\end{aligned}\end{align}
where $\iprod{\cdot,\cdot}_Y$ is the $L^2-$scalar product over $Y$. Set
\begin{align}\begin{aligned}\label{e.295}
V_\pm:=\set{\phi\in \Ker D_Y\,|\, \GS \,\phi=\pm\phi},
\end{aligned}\end{align}
which are two Lagrangian subspaces of $(\Ker D_Y,\Phi)$.
Let ${\rm pr}_{V_+}$ be the orthogonal projection from $L^2(Y,S|_Y)$ onto $V_+$. Set
\begin{align}\begin{aligned}\label{e.296}
\Pi_{V_+}:=P_>+\pr{V_+}.
\end{aligned}\end{align}
We impose the generalized Atiyah-Patodi-Singer global spectral projection boundary condition on $D_M$,
\begin{align}\begin{aligned}\label{e.297}
\Dom (D_{M,V_+}):=\set{\psi\in \Cinf{M,S}: \Pi_{V_+}\left(\psi|_Y\right)=0}.
\end{aligned}\end{align}
Then the operator $D_{M,V_+}:L^2(M,S)\ra L^2(M,S)$ is essentially self-adjoint (cf. \cite[Lemma 1.11]{Muller94}).

With respect to the $\Z_2-$grading over $S$, we can write the odd operator $D_{M,V_+}$ in matrix form,
\begin{align}\begin{aligned}\label{e.298}
D_{M,V_+}=\left(
      \begin{array}{cc}
        0 & D_{M,V_+}^-\\
        D_{M,V_+}^+ & 0 \\
      \end{array}
    \right).
\end{aligned}\end{align}
Since $D_{M,V_+}$ is self-adjoint, one has $\left(D^+_{M,V_+}\right)^*=D^-_{M,V_+}$. The generalized Atiyah-Patodi-Singer index of $D_{M,V_+}$ is defined as:
\begin{align}\begin{aligned}\label{e.299}
\ind(D_{M,V_+}):=\dim \Ker D^+_{M,V_+}-\dim \Ker D^-_{M,V_+}=\tr\left(\Gamma_S|_{\Ker (D_{M,V_+})}\right).
\end{aligned}\end{align}

\begin{figure}
  \centering
  \includegraphics[width=10cm]{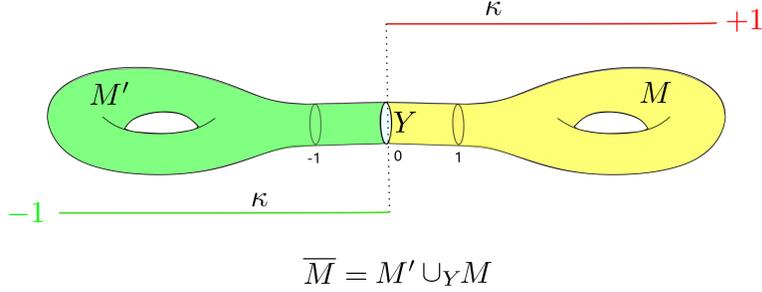}\\
\caption{Double manifold and domain-wall function}\label{figure11}
\end{figure}

Let $M'$ be another copy of $M$, such that $M'$ has $Y\times[-1,0]$ as product neighborhood of $Y=\partial M'$. We identify both the boundaries of $M,\, M'$ with $Y\times \set{0}$.
We construct the double manifold $\bM$ associated with $(M,Y)$ by gluing $M'$ and $M$ together along $Y\times \set{0}$ (see Figure \ref{figure11}), i.e.,
\begin{align}\begin{aligned}\label{e.300}
\bM=M'\cup_{Y}M.
\end{aligned}\end{align}
The domain-wall function $\kappa$ is a step function on the double manifold $\bM$, such that
\begin{align}\begin{aligned}\label{e.301}
\kappa(x)\equiv -1,\quad \text{for}\quad x\in M'\backslash Y;\quad \kappa(x)\equiv 1,\quad \text{for}\quad x\in M\backslash Y.
\end{aligned}\end{align}

By the product structures \eqref{e.17}, \eqref{e.18}, all the geometric objects on $M$, like $\gTM,\,h^S$, can be extended to the double manifold $\bM$. Let $D_{\bM}:\Cinf{\bM, S}\ra \Cinf{\bM,S}$ be the Dirac type operator over $\bM$, which is an extension of $D_M$. The operator $D_{\bM}$ is also an elliptic self-adjoint first order differential operator. Fix $m>0$, the massive Dirac-type operator is defined as
\begin{align}\begin{aligned}\label{e.302}
D_{\bM}-m\GS,
\end{aligned}\end{align}
which is also self-adjoint and elliptic, but no longer of odd parity. Let $\set{0\neq \lambda_j,j\in \N}$ be the set of eigenvalues of $D_{\bM}-m\GS$, then the eta function of $D_{\bM}-m\GS$ is defined to be
\begin{align}\begin{aligned}\label{e.303}
\eta(s,D_{\bM}-m\GS):=\sum_{\lj\neq 0}\frac{\sign(\lj)}{\abs{\lj}^s},\quad \Re(s)>m.
\end{aligned}\end{align}
It is absolutely convergent on $\set{s\in \C|\,\Re(s)>m}$ and can be extended to a meromorphic function on the whole complex plane \cite[(3.9)]{APS}.

The eta invariant associated to $D_{\bM}-m\GS$ is defined to be $\eta(D_{\bM}-m\GS):=\eta(0,D_{\bM}-m\GS)$. Similarly, we define the domain-wall Dirac type operator as
\begin{align}\begin{aligned}\label{e.304}
D_{\bM}+m\kappa\GS,
\end{aligned}\end{align}
and let $\eta(D_{\bM}+m\kappa\GS)$ be the eta invariant associated with $D_{\bM}+m\kappa\GS$ (see Definition \ref{d.2}).

The following theorem is one of the main results in this paper.

\begin{thm}\label{t.3}For $m>0$ sufficiently large and the Lagrangian subspace $V_+$ of $\Ker D_Y$ in \eqref{e.295}, the following identity holds
\begin{align}\begin{aligned}\label{e.47}
\ind(D_{M,V_+})=\frac{\eta(D_{\bM}+ m\kappa\Gamma_S)-\eta(D_{\bM}- m\Gamma_S)}{2}-n_+,
\end{aligned}\end{align}
where $n_+=\dim V_+$.
\end{thm}

Set
\begin{align}\begin{aligned}\label{e.305}
P_0:=\Id -P_<-P_>,\quad \text{and} \quad P_{\geq}:=P_>+P_0.
\end{aligned}\end{align}
The global spectral projection boundary condition originally used in \cite{APS} by Atiyah, Patodi and Singer is the following one,
\begin{align}\begin{aligned}\label{e.306}
\Dom (D_{M,P_{\geq}}):=\set{\psi\in \Cinf{M,S}: P_{\geq}\left(\psi|_Y\right)=0}.
\end{aligned}\end{align}
The obtained operator $D_{M,P_{\geq}}:L^2(M,S)\ra L^2(M,S)$ is elliptic, but not self-adjoint.
By using Proposition 21.4 in \cite{BW93}, we get the following theorem as a consequence of Theorem \ref{t.3}.

\begin{thm}\label{t.5}For $m>0$ sufficiently large and the Lagrangian subspace $V_+$ of $\Ker D_Y$ in \eqref{e.295}, the following identity holds
\begin{align}\begin{aligned}\label{e.181}
\ind(D_{M,P_{\geq}})=\frac{\eta(D_{\bM}+ m\kappa\Gamma_S)-\eta(D_{\bM}- m\Gamma_S)}{2}-2n_+,
\end{aligned}\end{align}
where $n_+=\dim V_+=\onehalf \dim \Ker D_Y$.
\end{thm}

Under the assumption that the boundary Dirac operator $D_Y$ is invertible, i.e., $\Ker D_Y=0$, Theorem \ref{t.5} was first proved in \cite{FF20}. In this paper we establish Theorem \ref{t.5} without the invertibility of
the boundary Dirac operator $D_Y$, and there is a new term appearing in the formula. This new term represents the contribution of the kernel of the Dirac operator $D_Y$ on the boundary.

Set
\begin{align}\begin{aligned}\label{e.307}
\widetilde{M}'=M'\backslash \Big(Y\times[-\onehalf,0]\Big),\quad,N=Y\times [-\onehalf,\onehalf],\quad \widetilde{M}=M\backslash \Big(Y\times [0,\onehalf]\Big).
\end{aligned}\end{align}
Then we decompose the double manifold $\bM$ into three pieces: $\widetilde{M}'$,\, $N$, \,$\widetilde{M}$, i.e.,
\begin{align}\begin{aligned}\label{e.308}
\bM=\widetilde{M}'\cup_Y N\cup_Y \widetilde{M}.
\end{aligned}\end{align}
Our strategy for proving Theorem \ref{t.3} is to break the difference
\begin{align}\begin{aligned}\label{e.309}
\eta(D_{\bM}+ m\kappa\Gamma_S)-\eta(D_{\bM}- m\Gamma_S)
\end{aligned}\end{align}
down into three pieces corresponding to $\widetilde{M}'$,\, $N$, \,$\widetilde{M}$ respectively, by using the asymptotic gluing formula for eta invariants. The two eta invariants in the piece over $\widetilde{M}'$ cancel out, hence the contribution from $\widetilde{M}'$ is zero. The piece over $\widetilde{M}$ gives the generalized Atiyah-Patodi-Singer index on the left hand side of \eqref{e.181}. The piece over $N$ gives the term in \eqref{e.181} representing the contribution of $\Ker D_Y$.

When we apply the adiabatic limit to prove the asymptotic gluing formula for eta invariants, we need to vary the length
of the cylinder part near the separating hyper-surface $Y$. It is well known that the Fredholm index of geometric elliptic operators over smooth Riemannian manifolds
is rather stable under the perturbation of the metrics. But in general the eta invariant is sensitive to the change of metrics, and may have integer jumps if it happens that some eigenvalues cross the zero point. As we will see, the difference of massive eta invariants \eqref{e.308} becomes very stable when the length of cylinder part near the boundary varies, if the mass parameter $m>0$ is large enough. This partially explains why we could have such a formula relating the generalized Atiyah-Patodi-Singer index to the eta invariants in large mass limits.

This paper is organized as follows. In Section \ref{s1}, we recall the definition of Dirac type operators and eta invariants. In Section \ref{s2}, we introduce the generalized Atiyah-Patodi-Singer index and the eta invariant associated with the domain-wall massive Dirac operators. In Section \ref{s3}, we prove the asymptotic gluing formula for eta invariants by using adiabatic limit method.
In Section \ref{s4}, we prove Theorem \ref{t.3} and Theorem \ref{t.5} by using the asymptotic gluing formula.
\newline

\textbf{Acknowledgments.}
The author is indebted to Zhi Hu for useful conversations on the topics of massive Dirac operators.

\section{Massive Dirac operator and eta invariants}\label{s1}

In this section, we will recall the definition of some basic objects used in our paper, like Dirac type operators, eta invariants and massive Dirac type operators.

\subsection{Massive Dirac operators}
Let $M$ be an oriented compact smooth Riemannian manifold of even dimension $n$ (with or without boundary). Let $S$ be a complex vector bundle over $M$. Let $\GME$ be the space of smooth sections of $S$ over $M$.

Let $\gTM$ be the Riemannian metric on the tangent bundle, which induces a metric $g^{T^*M}$ on the cotangent bundle $T^*M\rac{\pi}M$. Let $h^S$ be a Hermitian metric on $S$.

We recall the definition of the principal symbol of a differential operator $P:\Cinf{M,S}\ra \Cinf{M,S}$ of order $k\geq 1$. If we have
\begin{align}\begin{aligned}\label{e.182}
P=\sum_{0\leq \abs{\alpha}\leq k}a_\alpha(x)\frac{\partial^{\abs{\alpha}}}{\partial x^\alpha}\quad
\text{and}\quad \xi =\sum _{\mu=1}^n\xi_\mu dx_\mu\in T^*X_x
\end{aligned}\end{align}
in local coordinates, then the principal symbol $\sigma_P\in \Cinf{T^*M,\pi^*\End(S)}$ is locally given by
\begin{align}\begin{aligned}\label{e.183}
\sigma_P(x,\xi):=i^k\sum_{\abs{\alpha}=k}a_\alpha(x)\xi^{\alpha}:S_x\ra S_x.
\end{aligned}\end{align}

\begin{defn}\label{d.3}
A first order differential operator $D_M: \Cinf{M,S}\ra \Cinf{M,S}$ is called \emph{an operator of Dirac type}, if the principal symbol of $D^2_M$ is given by
\begin{align}\begin{aligned}\label{e.5}
\sigma_{D_M^2}(x,\xi)=\sum^n_{i,j=1}g^{ij}(x)\xi_i\xi_j,
\end{aligned}\end{align}
where $g^{ij}(x)$ is the Riemannian metric tensor in local coordinates
 (cf. \cite[Def. 2.2]{BGV92}, \cite[\S 3]{BW93}). The operator $D^2_M$ is called the \emph{Dirac Laplacian}.
\end{defn}

Let $\CLM$ be the Clifford algebra bundle associated to $(TM, \gTM)$, whose fiber at $x\in M$ is the Clifford algebra associated to the Euclidean space $(T_xM, g^{T_xM})$.
Assume that the complex vector bundle $S$ carries the structure of a complex Clifford module bundle over $M$. The Clifford multiplication is a map
\begin{align}\begin{aligned}\label{e.184}
\bc: \Cinf{M,TM\otimes S}\ra \Cinf{M,S}.
\end{aligned}\end{align}
Let $\nS:\Cinf{M,S}\ra \Cinf{M,T^*M\otimes S}$ be a connection on $S$ and let $J:\Cinf{M,T^*M\otimes S}\cong \Cinf{M,TM\otimes S}$ denote the isomorphism of vector and covector fields.

 \begin{defn}\label{d.4}
 \begin{enumerate}
   \item A \emph{generalized Dirac operator} $D_M:\Cinf{X,S}\ra \Cinf{X,S}$ is a first order differential operator over $S$ defined by
\begin{align}\begin{aligned}\label{e.185}
D=\bc\circ J\circ \nS.
\end{aligned}\end{align}
In terms of an orthonormal basis $\set{e_i|\,1\leq i\leq n}$ of $T_xM$, we have for $s\in \Cinf{M,S}$
\begin{align}\begin{aligned}\label{e.186}
Ds(x)=\sum^n_{i=1}\bc(e_i)\left(\nS_{e_i}s\right)_x.
\end{aligned}\end{align}
   \item If the connection $\nS$ is compatible with the $\CLM-$module structure of $S$ and the induced Levi-Civita connection $\n^{TM}$ on $\CLM$, i.e.,
       \begin{align}\begin{aligned}\label{e.187}
[\nS,\bc(e_i)]=\bc(\n^{TM}e_i), \quad \text{for}\quad 1\leq i\leq n,
\end{aligned}\end{align}
        we call $D_M$ a \emph{ compatible Dirac operator}.
 \end{enumerate}
\end{defn}

Let $S=S_+\oplus S_-$ be a $\Z_2-$graded structure over $S$. Let $\Gamma_S$ be the chirality operator on $S=\Sp\oplus\Sm$ such that $\Gamma_S|_{S_{\pm}}=\pm\Id$. We say that the Dirac-type operator is of odd parity, if
\begin{align}\begin{aligned}\label{e.189}
\GS D_M=-D_M \GS.
\end{aligned}\end{align}
 Assume that $h^S$ is a Hermitian metric on $S$ such that $S_+$ and $S_-$ are orthogonal. Note that if $S$ has a Clifford-module structure, there is a nature chirality operator
\begin{align}\begin{aligned}\label{e.188}
\GS=i^p\bc(e_1)\cdot \bc(e_2)\cdots \bc(e_n),
\end{aligned}\end{align}
where $p=n/2$ if $n$ is even, $p=(n+1)/2$ if $n$ is odd, and $\set{e_i,\,1\leq i\leq n}$ is an oriented orthonormal local frame of $TX$. One can verify $\GS^2=1$ for the chirality operator in \eqref{e.188}, which induces a natural $\Z_2$-splitting on $S$.

\begin{defn}
 Let $D_M$ be a Dirac-type operator of odd parity. For $m>0$, the massive Dirac operators are defined as
\begin{align}\begin{aligned}\label{e.7}
D_M\pm m\Gamma_S,
\end{aligned}\end{align}
which is no longer of odd parity.
\end{defn}

Since $\Gamma_S^2=1$ and $\Gamma_S D_M=-D_M\Gamma_S$, we have
\begin{align}\begin{aligned}\label{e.8}
(D_M\pm m\Gamma_S)^2=D_M^2+m^2.
\end{aligned}\end{align}
So the massive Dirac operator is still an operator of Dirac-type, and the associated Dirac Laplacian is a shift of the Dirac Laplacian $D_M^2$ by the square of mass $m$.

\subsection{Eta invariants on manifolds without boundary}

The Dirac operator $D_M:\GME\ra \GME$ is a first order elliptic differential operator on $M$, which is formally selfadjoint with respect to the $L^2$-product on $\GME$. Hence $D_M$ is essentially selfadjoint in $\LtME$. The eta invariant introduced by Atiyah, Patodi and Singer \cite{APS} is a non-local spectral invariant associated with $D_M$.

Let $\set{\lambda_j,j\in \N}$ be the set of eigenvalues of $D_M$, then the eta function of $D_M$ is defined to be
\begin{align}\begin{aligned}\label{e.9}
\eta(s,D_M):=\sum_{\lj\neq 0}\frac{\sign(\lj)}{\abs{\lj}^s},\quad \Re(s)>m.
\end{aligned}\end{align}
It is absolutely convergent on $\set{s\in \C|\,\Re(s)>m}$ and can be extended to a meromorphic function on the whole complex plane. By using the Mellin transform (cf. \cite[\S 9.6]{BGV92}), we can express the eta function in terms of the trace of the heat operator,
\begin{align}\begin{aligned}\label{e.10}
\eta(s,D_M)=\frac{1}{\G((s+1)/2)}\int_0^\infty t^{(s-1)/2}\tr(D_Me^{-tD_M^2})dt.
\end{aligned}\end{align}
It is well known that $\eta(s,D_M)$ is holomorphic at $s=0$ (cf. \cite{APSIII}, \cite{Gil81}, \cite{BG92}). The eta invariant of Dirac operator $D_M$ is defined to be
\begin{align}\begin{aligned}\label{e.11}
\eta(D_M):=\eta(0,D_M).
\end{aligned}\end{align}
It measures the spectral asymmetry of $D_M$, and appears naturally in the Atiyah-Patodi-Singer index theorem as a boundary correction term (cf. \cite{APS}). For compatible Dirac operators, $\eta(s,D_{M})$ is regular in the half-plane ${\rm Re}(s)>-2$ (cf. \cite[Theorem 2.6]{BF-II}). By \eqref{e.10} and $\G(\frac{1}{2})=\sqrt{\pi}$, the eta invariant of $D_{M}$ is given by
\begin{align}\begin{aligned}\label{e.12}
\eta(D_M)=\frac{1}{\sqrt{\pi}}\int_0^\infty t^{-1/2}\tr(D_Me^{-tD_M^2})dt,
\end{aligned}\end{align}
when $D_{M}$ is a compatible Dirac operator.

\subsection{Index and Eta invariants}
Assume that $M$ has no boundary. The Atiyah-Singer index of $D_M$ is defined as the Fredholm index of $D^+_M:=D_M|_{\Cinf{M,S_+}}$,
\begin{align}\begin{aligned}\label{e.13}
\ind(D_M):
&=\dim \Ker (D_M^+)-\dim {\rm Coker} (D_M^+)\\
&=\dim \Ker (D_M^+)-\dim \Ker (D_M^-).
\end{aligned}\end{align}
The famous Mckean-Singer formula says (cf. \cite{MS67}, \cite[Thm. 3.50]{BGV92}) that for any $t>0$
\begin{align}\begin{aligned}\label{e.14}
\ind(D_M)=\str(e^{-tD_M^2})=\tr(\Gamma_S e^{-tD_M^2})=\int_M\str(e^{-tD_M^2}(x,x))dv_x,
\end{aligned}\end{align}
where $dv$ is the Riemannian volume form. Here $\str$ denotes the super-trace, defined by
\begin{align}\begin{aligned}\label{e.190}
\str(\cdot):=\tr(\GS\,\cdot).
\end{aligned}\end{align}

Let $D_M\pm m\Gamma_S$ be the massive Dirac operators defined in \eqref{e.7} for $m>0$, which is not compatible. Then we have the following theorem.

\begin{thm}\label{t.1} The following identity holds
\begin{align}\begin{aligned}\label{e.15}
\ind(D_M)=\frac{\eta(D_M+m\GS)-\eta(D_M-m\GS)}{2}.
\end{aligned}\end{align}
\end{thm}
\begin{proof}
For general Dirac-type operators, not necessarily being compatible, by\eqref{e.8}, \eqref{e.10} and \eqref{e.14} we have
 \begin{align}\begin{aligned}\label{e.74}
&\qquad\qquad\onehalf(\eta(s,D_M+m\GS)-\eta(s,D_M-m\GS))\\
&=\frac{1}{2\G((s+1)/2)}\int_0^\infty t^{(s-1)/2}\left[\tr((D_M+m\GS)e^{-t(D_M^2+m^2)})\right.\\
&\qquad \qquad\qquad\qquad \qquad\qquad \qquad\qquad\qquad\left.-\tr((D_M-m\GS)e^{-t(D_M^2+m^2)})\right]dt\\
&=\frac{m}{\G((s+1)/2)}\int_0^\infty t^{(s-1)/2}\left[\tr(\GS e^{-t(D_M^2+m^2)})\right]dt\\
&=\frac{m\ind(D_M)}{\G((s+1)/2)}\int_0^\infty t^{(s-1)/2}e^{-tm^2}dt=\frac{\ind(D_M)}{m^{s}}.
\end{aligned}\end{align}
 Let $s\in \C$ go to zero in \eqref{e.74}, we get Equation \eqref{e.15}.
The proof of Theorem \ref{t.1} is completed.
\end{proof}

The main purpose of this paper is to generalize the above theorem to manifolds with boundaries. And we will prove a formula expressing the Atiyah-Patodi-Singer index in terms of eta invariants associated to massive Dirac operators and domain-wall massive Dirac operators.

\section{Domain-wall eta invariants}\label{s2}

In this section, we will introduce the generalized Atiyah-Patodi-Singer spectral projection boundary condition. Then we give the definition of generalized Atiyah-Patodi-Singer index and the definition of eta invariant associated to the domain-wall massive Dirac type operator.

\subsection{Atiyah-Patodi-Singer index}

Let $M$ be an even dimensional smooth manifold with smooth boundary $\partial M=Y$. Let $S=S_+\oplus S_-$ be a $\Z_2-$graded vector bundle over $M$. We denote the $\Z_2-$grading operator by $\Gamma_S$, i.e., $\Gamma_S |_{S_{\pm}}=\pm 1$.

 Let $\gTM$ be a Riemannian metric on $M$. Let $h^S$ be a Hermitian metric over $S$, such that $S_+$ and
$S_-$ are orthogonal.

\begin{defn}\label{d.1} For a compact manifold $X$ and a subset $I$ of $\R$, we set $X_I:=X\times I$, for example $X_{[-a,a]}=X\times [-a,a]$ ($a>0$), $X_{\R}=X\times (-\infty,+\infty)$.
\end{defn}

Recall that $\Cinf{M,S}$ denote the space of smooth sections of $S$. We denote $\Cinfc{M,S}$ the subspace of $\Cinf{M,S}$ consisting of sections with compact support in the interior of $M$. For $s,s'\in \Cinf{M,S}$, let
\begin{align}\begin{aligned}\label{e.36}
\iprod{s,s'}:=\int_M\iprod{s(x),s'(x)}_{h^S}dV_x
\end{aligned}\end{align}
be the $L^2-$inner product induced by $\gTM$ and $h^S$. We denote by $L^2(M,S)$ the $L^2-$completion of $\Cinfc{M,S}$.

The Dirac-type operator $D_M:\Cinf{M,S_\pm}\ra \Cinf{M,S_\mp}$ is a first-order elliptic differential operator of odd parity. We assume that $D_M$ is formally self-adjoint, i.e., for $s, s'\in \Cinfc{M,S}$
\begin{align}\begin{aligned}\label{e.37}
\iprod{D_Ms,s'}=\iprod{s,D_M s'}.
\end{aligned}\end{align}
Assume that on a collor neighborhood $[0,1)\times Y$ of the boundary $Y$ we have
\begin{align}\begin{aligned}\label{e.38}
D_M=\gamma\Big(\pd{u}{\,}+D_Y\Big),
\end{aligned}\end{align}
where $\gamma:S|_Y\ra S|_Y$ is a unitary bundle isomorphism and $D_Y:\Cinf{Y,S|_Y}\ra \Cinf{Y,S|_Y}$ is a first order elliptic operator on $Y$, such that
\begin{align}\begin{aligned}\label{e.39}
&\gamma^2=-\Id,\quad \gamma\Gamma_S=-\Gamma_S\gamma,\quad\gamma^*=-\gamma,\\
&\quad D_Y\gamma=-\gamma D_Y,\quad D_Y\Gamma_S=\Gamma_S D_Y.
\end{aligned}\end{align}
Let $D^*_Y$ be the formal adjoint of $D_Y$. We assume that $D_Y$ is self-adjoint, i.e., $D^*_Y=D_Y$.

Now we have a symmetric unbounded operator $D_M$ in $L^2(M,S)$ with domain $\Cinfc{M,S}$. We use the generalized Atiyah-Patodi-Singer global boundary conditions (cf. \cite{BW93}) to get a self-adjoint extension of $D_M$.

Let $P_>$ (resp. $P_<$) be the positive (resp. negative) spectral projection of $D_Y$, i.e., the orthogonal projection from $L^2(Y,S|_Y)$ onto the subspace spanned by all eigenvectors of positive (resp. negative) eigenvalues. Set
\begin{align}\begin{aligned}\label{e.196}
P_0:=\Id -P_<-P_>,\quad \text{and} \quad P_{\geq}:=P_>+P_0.
\end{aligned}\end{align}
It is well-known that the spectral projection $P_{\geq}$ is a pseudo-differential operator of order
zero (cf. \cite{APS}, \cite[Prop. 14.2]{BW93}).

 Set $V:=\Ker D_Y\subset L^2(Y,S|_Y)$. The isomorphism $\gamma$ induces a symplectic structure on $V$, that is
\begin{align}\begin{aligned}\label{e.40}
\Phi(x,y):=\iprod{\gamma x,y},
\end{aligned}\end{align}
where $\iprod{\cdot,\cdot}$ is the $L^2-$scalar product. A Lagrangian subspace $L$ of $V$ is a subspace satisfying
\begin{align}\begin{aligned}\label{e.41}
L\oplus \gamma L=V,\quad \Phi(L,L)=0.
\end{aligned}\end{align}
Let ${\rm pr}_L$ be the orthogonal projection from $L^2(Y,S|_Y)$ onto $L$. Set
\begin{align}\begin{aligned}\label{e.121}
\Pi_L:=P_>+\pr{L}.
\end{aligned}\end{align}
Since $\Ker D_Y$ consists of smooth sections, the projection $\pr{L}$ has a smooth kernel.
We define the domain of $D_M$ to be
\begin{align}\begin{aligned}\label{e.42}
\Dom (D_{M,L}):=\set{\psi\in \Cinf{M,S}: \Pi_L\left(\psi|_Y\right)=0}.
\end{aligned}\end{align}
Then $D_M$ has a self-adjoint extension to $L^2(M,S)$, which we will denote by $D_{M,L}$ (cf. \cite{DouWoj91}, \cite{LW96}, \cite{Muller94}).

Now we introduce a special Lagrangian subspace of $(V,\Phi)$, which plays an important role in our paper.
By \eqref{e.39}, we have
\begin{align}\begin{aligned}\label{e.73}
\GS D_Y=D_Y\GS,\quad \gamma\GS=-\GS\gamma.
\end{aligned}\end{align}
Hence $\GS$ induces an endomorphism on $V=\Ker D_Y$.
Let $V_\pm\subset V$ be the subspaces of $V=\Ker D_Y$ associated to $\pm 1$-eigenvalue of $\GS$ respectively, i.e.,
\begin{align}\begin{aligned}\label{e.70}
V_\pm:=\set{\phi\in V\,|\, \GS \,\phi=\pm\phi}.
\end{aligned}\end{align}
Moreover, by the second equation of \eqref{e.73} we have
\begin{align}\begin{aligned}\label{e.71}
\gamma(V_+)=V_-,\quad \gamma(V_-)=V_+.
\end{aligned}\end{align}
Since the sub-bundles $S_+$ and $S_-$ are orthogonal with respect to the Hermitian metric $h^S$ over $S$, the subspaces $V_+$ and $V_-$ are orthogonal with respect to the $L^2$-metric on $V$. Then $V_\pm$ are Lagrangian subspaces of $(V,\Phi)$ by \eqref{e.40}, \eqref{e.41} and \eqref{e.71}.
The Lagrangian subspaces $V_{\pm}$ play a special role in the proof of Theorem \ref{t.3}.

Let $\Pi_{V_+}$ be the projection operator in \eqref{e.121}. With respect to the $\Z_2-$splitting $S=S_+\oplus S_-$, we can write the odd operator $D_{M,V_+}$ in matrix form, that
is,
\begin{align}\begin{aligned}\label{e.43}
D_{M,V_+}=\left(
      \begin{array}{cc}
        0 & D_{M,V_+}^-\\
        D_{M,V_+}^+ & 0 \\
      \end{array}
    \right).
\end{aligned}\end{align}
Since $D_{M,V_+}$ is self-adjoint, one has $\left(D^+_{M,V_+}\right)^*=D^-_{M,V_+}$. The index of $D_{M,V_+}$ is defined as:
\begin{align}\begin{aligned}\label{e.44}
\ind(D_{M,V_+}):=\dim \Ker D^+_{M,V_+}-\dim \Ker D^-_{M,V_+}=\tr\left(\Gamma_S|_{\Ker (D_{M,V_+})}\right).
\end{aligned}\end{align}
By the Mckean-Singer formula, we have for $t>0$
\begin{align}\begin{aligned}\label{e.45}
\ind(D_{M,V_+})=\tr\left(\Gamma_S \cdot e^{-t D^2_{M,V_+}}\right)=\int_M \tr\left(\Gamma_S \cdot e^{-t D^2_{M,V_+}}(x,x)\right)dV_x.
\end{aligned}\end{align}
\newline

By \eqref{e.73} and \eqref{e.70}, we see that
\begin{align}\begin{aligned}\label{e.271}
\Pi_{V_+}\cdot\GS=\GS\cdot\Pi_{V_+}.
\end{aligned}\end{align}
By \eqref{e.271}, the operator $\GS$ preserves the domain of $D_{M,V_+}$, hence we have
\begin{align}\begin{aligned}\label{e.272}
(D_{M}\pm m\GS)_{V_+}&=D_{M,V_+}\pm m\GS,\\
(D_{M}\pm m\GS)_{V_+}^2&=D_{M,V_+}^2\pm mD_{M,V_+}\circ\GS\pm m\GS \circ D_{M,V_+}+m^2\GS^2\\
&=D_{M,V_+}^2+m^2.
\end{aligned}\end{align}
Note that \eqref{e.272} does not hold for general Lagrangian subspaces of $(V,\Phi)$.
For $m>0$, the eta function associated to $D_{M,V_+}\pm m\GS$ is given by
\begin{align}\begin{aligned}\label{e.122}
&\eta(s,D_{M,V_+}\pm m\GS)\\
=& \frac{1}{\G((s+1)/2)}\int_0^\infty t^{(s-1)/2}\tr\left[(D_{M,V_+}\pm m\GS)e^{-t(D_{M,V_+}\pm m\GS)^2}\right]dt\\
=&\frac{1}{\G((s+1)/2)}\int_0^\infty t^{(s-1)/2}e^{-m^2t}\tr\left(D_{M,V_+}e^{-tD_{M,V_+}^2}\right)dt\pm \frac{\ind(D_{M,V_+})}{m^{s}}&,
\end{aligned}\end{align}
where we have used \eqref{e.45} in the last equality. Similarly to Theorem \ref{t.1}, we have the following theorem.

\begin{thm}\label{t.2}
For $m>0$ and the Lagrangian subspace $V_+$ of $V$, the following identity holds
\begin{align}\begin{aligned}\label{e.46}
\ind(D_{M,V_+})=\frac{\eta(D_{M,V_+}+ m\GS)-\eta(D_{M,V_+}- m\GS)}{2}.
\end{aligned}\end{align}
\end{thm}
\begin{proof}
Equation \eqref{e.46} follows from \eqref{e.8}, \eqref{e.10}, \eqref{e.45} exactly by the same arguments of Theorem \ref{t.1}, since by \eqref{e.122} we have
\begin{align}\begin{aligned}\label{e.273}
\onehalf(\eta(s,D_{M,V_+}+ m\GS)-\eta(s,D_{M,V_+}- m\GS))=\frac{\ind(D_{M,V_+})}{m^{s}}.
\end{aligned}\end{align}
Let $s$ go to zero, we get \eqref{e.46}. The proof of Theorem \ref{t.2} is completed.
\end{proof}

\subsection{Computations on the cylinders}
On the product $Y_{R^+}$ of $Y$ with half-line $u\geq 0$, for $m>0$ we consider the operator
\begin{align}\begin{aligned}\label{e.123}
D_{Y_{R^+}}-m\GS=\gamma\Big(\partial_u+D_Y+m\gamma\GS\Big).
\end{aligned}\end{align}
 By \eqref{e.70} and \eqref{e.71}, we could introduce two spectral projection operators as
\begin{align}\begin{aligned}\label{e.124}
\Pi_{V_+}=P_>+\pr{V_+},\qquad\Pi_{V_-}=P_<+\pr{V_-}.
\end{aligned}\end{align}
By \eqref{e.39}, \eqref{e.73}, \eqref{e.271}, we consider
\begin{align}\begin{aligned}\label{e.125}
(D_{Y_{R^+}}-m\GS)_{V_+}^2=(D_{Y_{R^+},V_+})^2+m^2=-\partial^2_u+D^2_Y+m^2,
\end{aligned}\end{align}
with the boundary condition
\begin{align}\begin{aligned}\label{e.126}
\Pi_{V_+}f(\cdot,0)=0 \quad \text{and}\quad \Pi_{V_-}\left\{\left(\pd{u}{f}+D_Yf+m\gamma\GS f\right)_{u=0}\right\}=0.
\end{aligned}\end{align}

Let $\set{\phi_{\lambda,i}|\,1 \leq i\leq \dim E_\lambda,\,\lambda>0}$ be an orthonormal basis of the space $E_\lambda$ of eigensections of $D_Y$ associated to the positive eigenvalue $\lambda>0$, such that
 \begin{align}\begin{aligned}\label{e.131}
D_Y\phi_{\lambda,i}=\lambda \phi_{\lambda,i}.
\end{aligned}\end{align}
Then the set $\set{\gamma\phi_{\lambda,i},\,\lambda>0}$ forms an orthonormal basis of the space $E_{-\lambda}$ for $\lambda>0$ by the second and fourth equations of \eqref{e.39}, such that
 \begin{align}\begin{aligned}\label{e.132}
D_Y\gamma\phi_{\lambda,i}=-\lambda \gamma\phi_{\lambda,i}.
\end{aligned}\end{align}
Set $n_+=\dim V_+$. Let $\set{\psi_i, 1\leq i\leq n_+}$ be an orthonormal basis of $V_+$, then the set
 \begin{align}\begin{aligned}\label{e.127}
\set{\phi_{\lambda_k},\,\gamma\phi_{\lambda_k},\,\psi_i,\,\gamma\psi_i|\,0<\lambda_k\in {\rm Spec}(D_Y),\,
k\in \N,\,1\leq i\leq n_+}
\end{aligned}\end{align}
gives an orthonormal basis of $L^2(Y,S|_Y)$.
\newline

By \cite[(2.16),(2.17)]{APS}, \cite[(1.13)]{Muller94}, the kernel $e^{-t(D_{Y_{R^+}}-m\GS)_{V_+}^2}((u,x),(v,y))$ of
$e^{-t(D_{Y_{R^+}}-m\GS)_{V_+}^2}=e^{-tm^2}e^{-t(D_{Y_{R^+},V_+})^2}$
 with the boundary condition \eqref{e.126} is given by
\begin{align}\begin{aligned}\label{e.142}
&e^{-m^2t}\sum_{\lambda_k>0}\frac{e^{-\lambda_k^2t}}{\sqrt{4\pi t}}\left\{\exp\left(\frac{-(u-v)^2}{4t}\right)-\exp\left(\frac{-(u+v)^2}{4t}\right)\right\}
\phi_{\lambda_k}(x)\otimes \overline{\phi_{\lambda_k}}(y)\\
&+e^{-m^2t}\sum_{\lambda_k>0}\left\{\frac{e^{-\lambda_k^2t}}{\sqrt{4\pi t}}\left(\exp\left(\frac{-(u-v)^2}{4t}\right)+\exp\left(\frac{-(u+v)^2}{4t}\right)\right)\right.\\
&\left.\qquad\qquad\qquad\qquad-\lambda_k e^{\lambda_k(u+v)}\erfc\left(\frac{u+v}{2\sqrt{t}}+\lambda_k
\sqrt{t}\right)\right\}\gamma\phi_{\lambda_k}(x)\otimes \overline{\gamma\phi_{\lambda_k}}(y)
\\
&+e^{-m^2t}\frac{1}{\sqrt{4\pi t}}\left\{\exp\left(\frac{-(u-v)^2}{4t}\right)-\exp\left(\frac{-(u+v)^2}{4t}\right)\right\}
\sum_{i=1}^{n_+}\psi_i(x)\otimes \overline{\psi_i}(y)\\
&+e^{-m^2t}\frac{1}{\sqrt{4\pi t}}\left\{\exp\left(\frac{-(u-v)^2}{4t}\right)+\exp\left(\frac{-(u+v)^2}{4t}\right)\right\}\sum_{i=1}^{n_+}\gamma\psi_i(x)\otimes \overline{\gamma\psi_i}(y),
\end{aligned}\end{align}
where $\set{\phi_{\lambda_i},\, i\in \N}$ is an orthonormal basis of ${\rm Ran}(P_>)$ consisting of the eigensections of $D_Y$ with eigenvalues $0<\lambda_1\leq \lambda_2\leq \cdots$, $\set{\psi_i,\,1\leq i\leq n_+}$ is an orthonormal basis of $V_+\subset \Ker D_Y$ and
$\erfc(x)$ is the complementary error function defined by
\begin{align}\begin{aligned}\label{e.144}
\erfc(x)=\frac{2}{\sqrt{\pi}}\int_x^\infty e^{-u^2}du.
\end{aligned}\end{align}
\newline

Now we consider the same problem for $D_{Y_{\R^-}}=\widetilde{\gamma}\left(\partial_u+\widetilde{D_Y}\right)$ on $Y_{\R^-}\cong (-\infty,0]\times Y$.
Note that $\widetilde{\gamma}$ is equal to $-\gamma$ in this case because of the orientation, and we have $\widetilde{D_Y}=-D_Y$, $\widetilde{\sigma}=-\sigma$. Now we consider the heat kernel of $e^{-t(D_{Y_{\R^-}}-m\GS)^2}$ with the boundary condition
\begin{align}\begin{aligned}\label{e.195}
\tPi_{V_-}f(\cdot,0)=0 \quad \text{and}\quad \tPi_{V_+}\left\{\left(\pd{u}{f}+\widetilde{D_Y}f+m\widetilde{\gamma }\GS f\right)_{u=0}\right\}=0,
\end{aligned}\end{align}
where $\tPi_{V_{\pm}}$ are the corresponding projection operators associated to $\widetilde{D_Y}$.
The generalized Atiyah-Patodi-Singer boundary condition \eqref{e.195} is equivalent to
\begin{align}\begin{aligned}\label{e.155}
\Pi_{V_+}f(\cdot,0)=0 \quad \text{and}\quad \Pi_{V_-}\left\{\left(\pd{u}{f}-D_Yf-m\gamma\GS f\right)_{u=0}\right\}=0.
\end{aligned}\end{align}

Similar to \eqref{e.142}, the kernel $e^{-t(D_{Y_{\R^-}}-m\GS)_{V_-}^2}((u,x),(v,y))$ of the heat operator $e^{-t(D_{Y_{\R^-}}-m\GS)_{V_-}^2}$ $=e^{-m^2t}e^{-t(D_{Y_{\R^-},V_-})^2}$ with the boundary condition \eqref{e.155} is given by
\begin{align}\begin{aligned}\label{e.160}
&e^{-m^2t}\sum_{\lambda_k>0}\frac{e^{-\lambda_k^2t}}{\sqrt{4\pi t}}\left\{\exp\left(\frac{-(u-v)^2}{4t}\right)-\exp\left(\frac{-(u+v)^2}{4t}\right)\right\}
\phi_{\lambda_k}(x)\otimes \overline{\phi_{\lambda_k}}(y)\\
&+e^{-m^2t}\sum_{\lambda_k>0}\left\{\frac{e^{-\lambda_k^2t}}{\sqrt{4\pi t}}\left(\exp\left(\frac{-(u-v)^2}{4t}\right)+\exp\left(\frac{-(u+v)^2}{4t}\right)\right)\right.\\
&\left.\qquad\qquad\qquad\qquad-\lambda_k e^{-\lambda_k(u+v)}\erfc\left(-\frac{u+v}{2\sqrt{t}}+\lambda_k
\sqrt{t}\right)\right\}\gamma\phi_{\lambda_k}(x)\otimes \overline{\gamma\phi_{\lambda_k}}(y)
\\
&+e^{-m^2t}\frac{1}{\sqrt{4\pi t}}\left\{\exp\left(\frac{-(u-v)^2}{4t}\right)-\exp\left(\frac{-(u+v)^2}{4t}\right)\right\}
\sum_{i=1}^{n_+}\psi_i(x)\otimes \overline{\psi_i}(y)\\
&+e^{-m^2t}\frac{1}{\sqrt{4\pi t}}\left\{\exp\left(\frac{-(u-v)^2}{4t}\right)+\exp\left(\frac{-(u+v)^2}{4t}\right)\right\}\sum_{i=1}^{n_+}\gamma\psi_i(x)\otimes \overline{\gamma\psi_i}(y),
\end{aligned}\end{align}
where $\set{\phi_{\lambda_i},\, i\in \N}$ is an orthonormal basis of ${\rm Ran}(P_>)$ consisting of the eigensections of $D_Y$ with eigenvalues $0<\lambda_1\leq \lambda_2\leq \cdots$, $\set{\psi_i,\,1\leq i\leq n_+}$ is an orthonormal basis of $V_+$.

\subsection{Domain-wall massive Dirac operator}
Let $M$ be an oriented even dimensional Riemannian manifold with a compact boundary $Y$. Let $m=\dim M$, then $Y$ is a smooth manifold of dimension $(m-1)$.

Let $D_Y$ be a Dirac-type operator on the boundary $Y$, which is a self-adjoint first-order elliptic operator, such that on the product neighborhood $\Yzo$ we have (see \eqref{e.38},\,\eqref{e.39})
\begin{align}\begin{aligned}\label{e.19}
D_M=\gamma\big(\frac{\partial\,}{\partial u}+D_Y\big).
\end{aligned}\end{align}

\begin{figure}
  \centering
  \includegraphics[width=10cm]{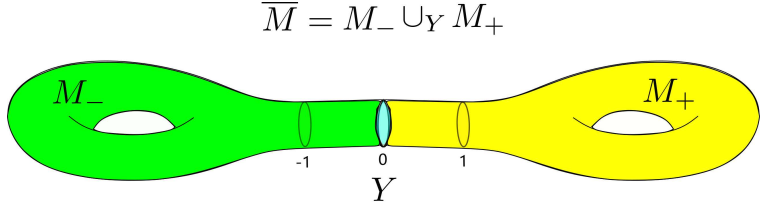}\\
\caption{Double manifold}\label{figure1}
\end{figure}

Let $\Mp,\Mm$ be two copies of $M$, such that $\Mp$ (resp. $\Mm$) has $\Yzo$ (resp. $\Ymo$) as product neighborhoods of $Y$. We identify both the boundaries $\partial M_{\pm}$ of $M_{\pm}$ with $Y_{\set{0}}$.
We construct the double manifold $\bM$ associated with $(M,Y)$ by gluing $\Mp$ and $\Mm$ together along $Y_{\set{0}}$ (see Figure \ref{figure1}), i.e.,
\begin{align}\begin{aligned}\label{e.20}
\bM=\Mm\cup_{Y}\Mp.
\end{aligned}\end{align}

Let $\kappa:[-1,1]\ra [-1,1]$ be a step function on $[-1,1]$, such that
\begin{align}\begin{aligned}\label{e.21}
\kappa(u)=\left\{
          \begin{array}{ll}
            -1, & \hbox{$u<0$;} \\
            1, & \hbox{$u>0$.}
          \end{array}
        \right.
\end{aligned}\end{align}
It can be naturally extended to $\Ymoo$, then a function on $\bM$, such that
$\kappa|_{M_{\pm}\backslash Y}=\pm 1$. We will call it the \emph{domain-wall function} on $\bM$.

\begin{defn}
  The \emph{domain-wall massive Dirac operator} on $\bM$ is defined as
\begin{align}\begin{aligned}\label{e.22}
\DDWm=D_{\bM}+m\kappa\Gamma_S,
\end{aligned}\end{align}
where $m>0$ and $D_{\bM}$ is the massless Dirac-type operator on $\bM$ (see \eqref{e.5}).
\end{defn}

For $m>0$, set
\begin{align}\begin{aligned}\label{e.23}
\DPVm=D_{\bM}-m\Gamma_S,
\end{aligned}\end{align}
in Physics language $\DPVm$ is called the \emph{Pauli-Villars massive Dirac operator} on $\bM$.

\subsection{Smooth approximation}

\begin{figure}
  \centering
  \includegraphics[width=6cm]{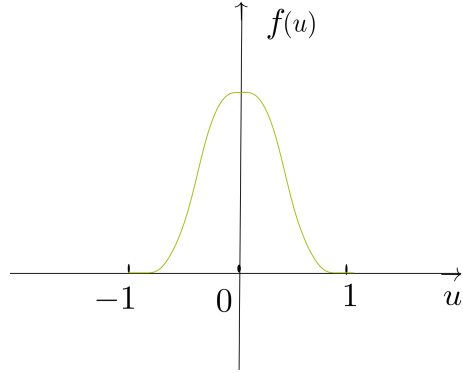}\\
\caption{Bump function with compact support in $(-1,1)$.}\label{figure2}
\end{figure}

\begin{figure}
  \centering
  \includegraphics[width=8cm]{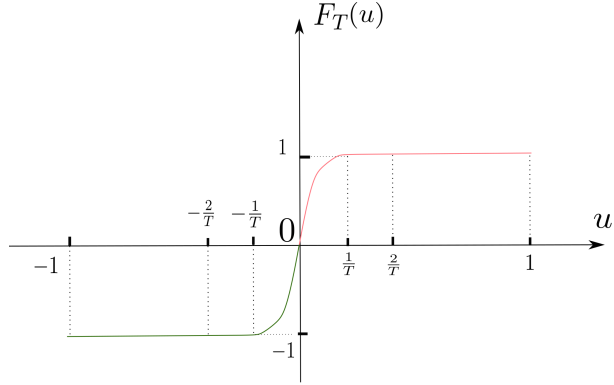}\\
\caption{Image of $F_T(u)$.}\label{figure5}
\end{figure}

Let $f(u),\,u\in (-\infty,\infty)$ be a bump function with compact support in $(-1,1)$, such that $f$ is even and
$\int_{-\infty}^\infty f(u)du=1$ (see Figure \ref{figure2}).
For $T>0$, we define
\begin{align}\begin{aligned}\label{e.26}
f_T(u)=Tf(Tu),
\end{aligned}\end{align}
then $f_T(u)$ is compactly supported in $(-\Tmo,\Tmo)$ and still we have $\int_{-\infty}^{\infty}f_T(u)du=1$. Set (see Figure \ref{figure5})
\begin{align}\begin{aligned}\label{e.27}
F_{T}(u)=2\int_{-\infty}^uf_T(t)dt-1.
\end{aligned}\end{align}
We extend $f_T(u),\, F_T(u)$ naturally to $\bM$. While $T\ra \infty$, we have
\begin{align}\begin{aligned}\label{e.28}
F_T(x)\ra \kappa(x).
\end{aligned}\end{align}

We define the deformed massive Dirac operator as
\begin{align}\begin{aligned}\label{e.29}
D_{\bM,m,T}=D_{\bM}+mF_T(x)\Gamma_S.
\end{aligned}\end{align}
The associated Laplacian is given by
\begin{align}\begin{aligned}\label{e.30}
D_{\bM,m,T}^2=D^2_{\bM}+2mf_T(x)\gamma\Gamma_S+m^2 F_T^2(x).
\end{aligned}\end{align}
And the eta function associated to $D_{\bM,m,T}$ is
\begin{align}\begin{aligned}\label{e.31}
\eta(s,D_{\bM,m,T})=\frac{1}{\G((s+1)/2)}\int_0^\infty t^{(s-1)/2}\tr(D_{\bM,m,T}e^{-tD_{\bM,m,T}^2})dt.
\end{aligned}\end{align}

\subsection{Domain-wall eta invariants}
Recall that $\bM$ is the double manifold of $M$.

 Since $\Ker(D_{\bM}+m\kappa \Gamma_S)=\set{0}$ for $m>0$ sufficiently large, there exists a constant $C_m>0$ such that
for any $f\in L^2(\bM;\End(S))$ with $\norm{f}_{L^2}<C_m$,
\begin{align}\begin{aligned}\label{e.33}
\Ker(D_{\bM}+m\kappa \Gamma_S+f)=\set{0}.
\end{aligned}\end{align}

\begin{lemma}
  Let $f_1\in L^2(\bM,\End(S))$ with $\norm{f_1}_{L^2}<C_m$ and $f_2\in L^2(\bM,\End(S))$ with $\norm{f_2}_{L^2}<C_m$. Assume that $m\kappa \Gamma_S+f_1$ and $m\kappa \Gamma_S+f_2$ are smooth operators. Then we have
  \begin{align}\begin{aligned}\label{e.35}
\eta(D_{\bM}+m\kappa \Gamma_S+f_1)=\eta(D_{\bM}+m\kappa \Gamma_S+f_2).
\end{aligned}\end{align}
\end{lemma}

\begin{defn}[Domain-wall $\eta$-invariant]\label{d.2}
  We define the eta invariant of domain-wall fermion Dirac operators by
  \begin{align}\begin{aligned}\label{e.34}
\eta(D_{\bM}+m\kappa \Gamma_S):=\eta(D_{\bM}+m\kappa\Gamma_S+f)
\end{aligned}\end{align}
for any $f\in L^2(\bM,\End(S))$ such that $\norm{f}_{L^2}<C_m$ and $m\kappa\Gamma_S+f$ is a smooth operator.
\end{defn}

 By Definition \ref{d.2} and \eqref{e.29}, for $T>0$ sufficiently large, we have
\begin{align}\begin{aligned}\label{e.32}
\eta(D_{\bM}+m\kappa \Gamma_S)=\eta(D_{\bM}+mF_T\Gamma_S).
\end{aligned}\end{align}

\section{An asymptotic gluing formula for eta invariants}\label{s3}
In this section we will prove a version of asymptotic gluing formula for eta invariants associated to massive Dirac-type operators. The gluing formula of eta invariants has been studied by many mathematicians (cf. \cite{bunke95}, \cite{Wojcie},\cite{Wo95}). Our proof is based on the splitting principle in the adiabatic limit developed by Douglas and Wojcieshowski (cf. \cite{DouWoj91}).

\subsection{Gluing problem and Adiabatic limit}

\begin{figure}
 \centering
  \includegraphics[width=8cm]{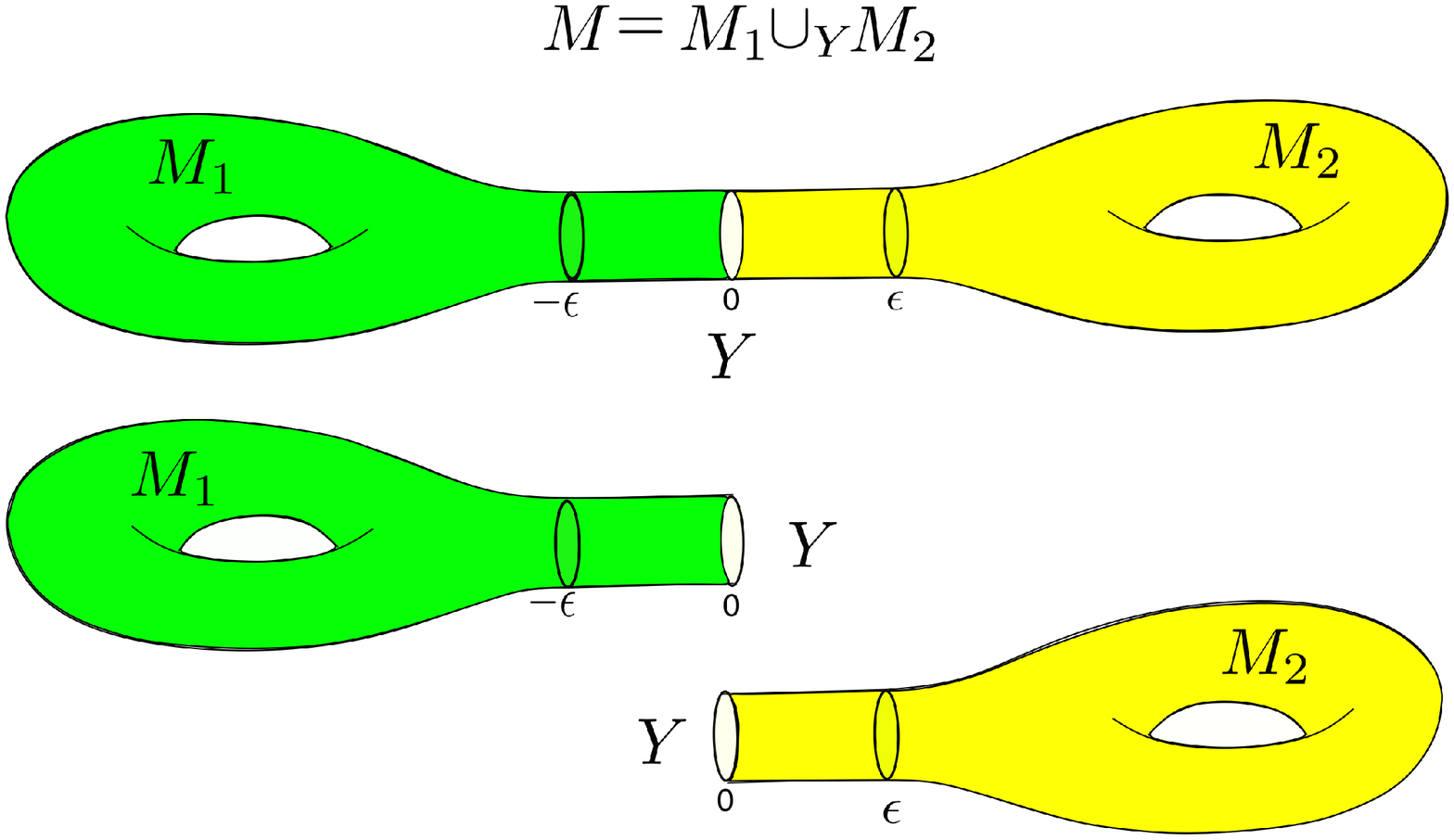}\\
\caption{Gluing setting}\label{figure7}
\end{figure}

Let $M$ be a compact Riemannian manifold without boundary. Assume that $M$ is split into two parts $M_1,\,M_2$ by a compact hypersurface $Y$, i.e., (see Figure \ref{figure7})
\begin{align}\begin{aligned}\label{e.75}
M=M_1\cup_Y M_2.
\end{aligned}\end{align}
Assume that $M$ has a cylinder part $Y_{[-\epsilon,\epsilon]}\cong Y\times [-\epsilon,\epsilon]$ for $\epsilon>0$, such that $\cY{-\epsilon,0}$ ( resp. $\cY{0,\epsilon}$ ) is the cylinder end of $M_1$ (resp. $M_2$).

Let $S$ be a complex vector bundle over $M$. Let $h^S$ denote the Hermitian metric on $S$. Assume that all metrics have product structures on the cylinder part $Y_{[-\epsilon,\epsilon]}$ (see \eqref{e.17}, \eqref{e.18}).

For $m>0$, let $D_{M,m}=D_{M}+mf$ be a family of Dirac-type operators, where $f$ is a smooth section of $\End(S)$. Assume that $D_{M,m}$ has the following form on the cylinder part $\cY{-\epsilon,\epsilon}$, that's
\begin{align}\begin{aligned}\label{e.76}
D_{M,m}|_{\cY{-\epsilon,\epsilon}}&=D_M-m\GS\\
&=\gamma\big(\frac{\partial\,}{\partial u}+D_Y+m\gamma\GS \big),
\end{aligned}\end{align}
where $D_M=\gamma\big(\frac{\partial\,}{\partial u}+D_Y\big)$ on $\cY{-\epsilon,\epsilon}$ (see \eqref{e.38}). We obtain two Dirac-type operators $D_{M_1},\,D_{M_2}$ by restricting $D_{M}$ to $M_1,\,M_2$ respectively.

Let $V:=\Ker(D_Y)$ with the symplectic structure in \eqref{e.40}.
Recall that $V_+,\,V_-$ are the Lagrangian subspaces of $V$ defined in \eqref{e.70}.
 Set
 \begin{align}\begin{aligned}\label{e.95}
 \Pi_{V_+}:=P_>+\pr{V_+} \quad \text{and}\quad \Pi_{V_-}:=P_<+\pr{V_-},
\end{aligned}\end{align}
where $P_<, \,P_>$ are spectral projection operators associated to $D_Y$ introduced in \eqref{e.196}.

  Let $\Domm,\,\Dtmp$ be the Dirac-type operators on $M_1,\, M_2$ with domains
\begin{align}\begin{aligned}\label{e.77}
&\Dom (\Domm):=\set{\psi\in \Cinf{M_1,S}: \Pi_{V_-}\left(\psi|_Y\right)=0},\\
&\Dom (\Dtmp):=\set{\psi\in \Cinf{M_2,S}: \Pi_{V_+}\left(\psi|_Y\right)=0}.
\end{aligned}\end{align}

Let $\eta(s,\Domm),\,\eta(s,\Dtmp)$ be the corresponding eta functions. The main concern of the gluing problem is to study
\begin{align}\begin{aligned}\label{e.78}
\delta_m(s;V_-,V_+):=\eta(s,D_{M}+mf)-\eta(s,\Domm)-\eta(s,\Dtmp).
\end{aligned}\end{align}
We will calculate the right hand side of \eqref{e.78} by the adiabatic limit method.

Let $\tM_1=M_1\backslash \cY{-\epsilon,0}$ and $\tM_2=M_2\backslash \cY{0,\epsilon}$. For $R>\epsilon>0$, we introduce some manifolds with stretched cylinder parts:
\begin{figure}
 \centering
  \includegraphics[width=8cm]{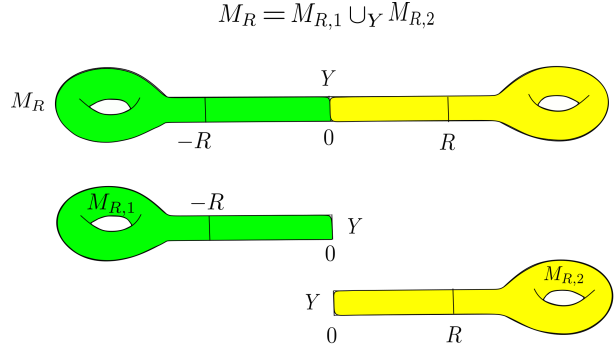}\\
\caption{Manifolds with stretched cylinder parts}\label{figure8}
\end{figure}

\begin{align}\begin{aligned}\label{e.79}
&\qquad M_R=\tM_1\cup_Y \cY{-R,R} \cup_Y \tM_2,\\
&M_{R,1}=\tM_1\cup_Y \cY{-R,0}, \qquad M_{R,2}= \cY{0,R}\cup_Y \tM_2.
\end{aligned}\end{align}
And we have (see Figure \ref{figure8})
\begin{align}\begin{aligned}\label{e.80}
& M_R=M_{R,1}\cup_Y M_{R,2}.
\end{aligned}\end{align}
By using the product structures in \eqref{e.17} and \eqref{e.18}, we can extend the complex vector bundle $S$ and the metrics naturally to these manifolds with stretched cylinder parts.

For ease of notations, we denote
\begin{align}\begin{aligned}\label{e.107}
&\qquad \qquad A_{R,m}:=D_{M_{R}}+mf,\\
&\quad A_{1,R,m}:=(D_{M_{R,1}}+mf)_{V_-},\quad A_{2,R,m}:=(D_{M_{R,2}}+mf)_{V_+}.
\end{aligned}\end{align}

Set
\begin{align}\begin{aligned}\label{e.81}
\delta_{m,R}(s;V_-,V_+):=\eta(s,A_{R,m})-\eta(s,A_{1,R,m})-\eta(s,A_{2,R,m}).
\end{aligned}\end{align}

Let
\begin{align}\begin{aligned}\label{e.89}
M_{\infty,1}=\tM_1\cup_Y \cY{0,\infty}\quad \text{and} \quad M_{\infty,2}= \cY{-\infty,0}\cup_Y \tM_2,
\end{aligned}\end{align}
be manifolds with infinite cylinder ends. The Dirac-type operators have natural extensions on these manifolds, which we denote by $D_{M_{\infty,1}},\,D_{M_{\infty,2}}$. As shown in \cite[Lemma 3.2]{Muller94}, the operators $D_{M_{\infty,1}},\,D_{M_{\infty,2}}$ are essentially self-adjoint. We denote by $\Ker_{L^2}(D_{M_{\infty,i}}+mf)$ the spaces of $L^2$-integrable solutions of the operators $D_{M_{\infty,i}}+mf$, for $i=1,2$.

\begin{Assumption}\label{a.1}
  We assume that for $m>0$ sufficiently large
\begin{align}\begin{aligned}\label{e.82}
\Ker_{L^2}(D_{M_{\infty,1}}+mf)=0,\quad \Ker_{L^2}(D_{M_{\infty,2}}+mf)=0,\quad \Ker(D_Y+m\gamma\GS)=0.
\end{aligned}\end{align}
\end{Assumption}

\begin{cor}\label{c.1}
  For $m>0$ sufficiently large, we have
  \begin{align}\begin{aligned}\label{e.88}
\Ker A_{R,m}=0,\quad \Ker A_{i,R,m}=0,
\end{aligned}\end{align}
where $i=1,2$ and $R\in [\epsilon,\infty)$.
\end{cor}
\begin{proof}
By \cite[Theorem 0.3]{Wojcie} and Assumption \ref{a.1}, there exists $R_0>0$ sufficiently large such that Equation \eqref{e.88} holds for any $R>R_0$. For $R\in [\epsilon,R_0]$, Equation \eqref{e.88} holds if we make $m>0$ sufficiently large.
\end{proof}

The variation of eta invariants is given by the integration of some local densities, which are locally computable from the jets of the symbol of the Dirac-type operators (cf. \cite[Prop. 2.8]{BF-II}, \cite{Gil93},\cite[Prop. 2.15]{Muller94}). Hence the difference of eta invariants $\delta_{m,R}(0;V_-,V_+)$ is independent of the parameter $R\in [\epsilon,\infty)$, which is half of the length of the cylinder part of $M_R$. Note that the eta invariant may have jumps in $2\Z$ while the eigenvalues cross the origin under the variation of the parameter $R\in [\epsilon,\infty)$, but Corollary \ref{c.1} implies that this phenomenon will not happen under the assumption \ref{a.1}.

\begin{prop}\label{p.1}
  Under the assumption \ref{a.1}, $\delta_{m,R}(0;V_-,V_+)$ does not depend on $R\in [\epsilon,\infty)$, i.e.,
\begin{align}\begin{aligned}\label{e.83}
\delta_m(0;V_-,V_+)=\delta_{m,R}(0;V_-,V_+),\quad \text{for any } R\in [\epsilon,\infty).
\end{aligned}\end{align}
\end{prop}

The same phenomenon as in Proposition \ref{p.1} has been described by Bunke \cite{bunke95} and Wojciechowski \cite{Wojcie}, in their proofs of the gluing formula for eta invariants in different situations. Generally, Equation \eqref{e.83} only holds modulo some integers. Since we are concerned with the index problem of Dirac-type operators, we need Assumption \ref{a.1} to make Equation \eqref{e.83} hold exactly in $\R$, not only in $\R/\Z$.

To prove the asymptotic gluing formula for massive Dirac-type operators, we only need to compute the limit of the right hand side of Equation \eqref{e.83} when $R$ goes to infinity. Douglas and Wojciechowski initiated the study of the adiabatic limit of eta invariants (cf. \cite{DouWoj91}). In fact, they developed a general method to calculate the adiabatic limits of the spectral invariants. We will follow their approach to establish the gluing formula of eta invariants for massive Dirac type operators.

We decompose the difference of eta functions $\delta_{m,R}(s;V_-,V_+)$ into the small time contribution and the large time contribution (cf. \cite{PW02}, \cite{PaWo06}, \cite{Wojcie}). Set
\begin{align}\begin{aligned}\label{e.84}
\Delta_{m,R}(t)=&\tr\left(A_{R,m}\exp(-tA_{R,m}^2)\right)-\sum_{i=1}^2\tr\left(A_{i,R,m}\exp(-t(A_{i,R,m}^2)\right).
\end{aligned}\end{align}
The small time contribution is defined as: for $\var>0$ small enough
\begin{align}\begin{aligned}\label{e.85}
\delta^S_{m,R}(s)=\frac{1}{\G((s+1)/2)}\int_0^{R^{2-\var}} t^{(s-1)/2}\Delta_{m,R}(t)dt.
\end{aligned}\end{align}
And the large time contribution is
\begin{align}\begin{aligned}\label{e.86}
\delta^L_{m,R}(s)=\frac{1}{\G((s+1)/2)}\int_{R^{2-\var}}^\infty t^{(s-1)/2}\Delta_{m,R}(t)dt.
\end{aligned}\end{align}

\subsection{Large time contributions}

First, we show that the large time contributions vanish when $R$ goes to infinity.

\begin{prop}\label{p.2}
Under Assumption \ref{a.1}, there exist $c_0>0,\,R_0>0$ such that for any eigenvalue $\mu_R$ of the operators $A_{R,m},\,A_{i,R,m},\,i=1,2$ and $R\in [R_0,\infty)$, we have
\begin{align}\begin{aligned}\label{e.90}
\abs{\mu_R}>c_0.
\end{aligned}\end{align}
\end{prop}
\begin{proof}
  Proposition \ref{p.2} follows from a similar argument as in \cite[Theorem 0.3]{Wojcie} under the assumption \ref{a.1}.
\end{proof}

\begin{thm}\label{t.4}
For $\var>0$ sufficiently small, the following equality holds
\begin{align}\begin{aligned}\label{e.87}
\lim_{R\ra 0}\delta^L_{m,R}(0)=0.
\end{aligned}\end{align}
\end{thm}

\begin{proof}
Let $\set{u_{k,R}}_{k=1}^\infty$ be the set of eigenvalues of $A_{R,m}$ for $R>R_0$.
\begin{align}\begin{aligned}\label{e.91}
&\frac{1}{\G((s+1)/2)}\int_{R^{2-\var}}^\infty t^{(s-1)/2}\tr\left(A_{R,m}\exp(-tA_{R,m}^2)\right)dt\\
=&\frac{1}{\G((s+1)/2)}\int_{R^{2-\var}}^\infty t^{(s-1)/2}\sum_{k=1}^\infty u_{k,R}\exp(-tu_{k,R}^2)dt.
\end{aligned}\end{align}
Since there exists a constant $C_1>0$ such that the function $xe^{-x^2}$ is uniformly bounded by $C_1$ for $x\in[0,\infty)$, the right hand side of \eqref{e.91} can be controlled by
\begin{align}\begin{aligned}\label{e.92}
&\frac{C_1}{\G((s+1)/2)}\int_{R^{2-\var}}^\infty t^{(s-1)/2}\sum_{k=1}^\infty \exp(-(t-1)u_{k,R}^2)dt\\
=&\frac{C_1}{\G((s+1)/2)}\int_{R^{2-\var}}^\infty t^{(s-1)/2}\sum_{k=1}^\infty \exp(-(t-2)u_{k,R}^2) \exp(-u_{k,R}^2)dt\\
\leq &\frac{C_1}{\G((s+1)/2)}\int_{R^{2-\var}}^\infty t^{(s-1)/2}\exp(-(t-2)c^2_0)\sum_{k=1}^\infty  \exp(-u_{k,R}^2)dt\\
\leq &\frac{C_2}{\G((s+1)/2)}\tr(\exp(-A_{R,m}^2))\int_{R^{2-\var}}^\infty t^{(s-1)/2}e^{-tc^2_0}dt\\
\leq &\frac{C_3 R}{\G((s+1)/2)}c^{-(s+1)}_0\int_{c^2_0R^{2-\var}}^\infty t^{(s-1)/2}e^{-t}dt.
\end{aligned}\end{align}
The same results hold for $A_{i,R,m},\,i=1,2$ by Proposition \ref{p.2}. Hence we have
\begin{align}\begin{aligned}\label{e.93}
\abs{\delta^L_{m,R}(s)}
\leq &\frac{C_4 R}{\G((s+1)/2)}c^{-(s+1)}_0\int_{c^2_0R^{2-\var}}^\infty t^{(s-1)/2}\exp(-t)dt.
\end{aligned}\end{align}
Put $s=0$ in \eqref{e.93}, we get
\begin{align}\begin{aligned}\label{e.94}
\abs{\delta^L_{m,R}(0)}
\leq &\frac{C_4 R}{c_0\sqrt{\pi}}\int_{c^2_0R^{2-\var}}^\infty t^{-1/2}\exp(-t)dt\\
\leq &\frac{C_4 R}{c_0\sqrt{\pi}}c_0^{-1}R^{-1+\var/2}\int_{c^2_0R^{2-\var}}^\infty \exp(-t)dt\\
\leq &C_5 R^{\var/2}\exp(-c^2_0R^{2-\var})\,\ra 0,\, \text{as}\, R\ra \infty.
\end{aligned}\end{align}
Equation \eqref{e.87} follows from \eqref{e.94}. The proof of Theorem \ref{t.4} is completed.
\end{proof}

\subsection{Small time contributions}

Now we use Duhamel's principle to handle the small time contribution $\delta^S_{m,R}(s)$.

Let $\rho(a,b): [0,\infty)\ra [0,1]$ be a cut-off function such that
\begin{align}\begin{aligned}\label{e.96}
\rho(a,b)(u)=\left\{
               \begin{array}{ll}
                 0, & \hbox{$0\leq u\leq a$;} \\
                 1, & \hbox{$u\geq b$.}
               \end{array}
             \right.
\end{aligned}\end{align}
Set
\begin{align}\begin{aligned}\label{e.97}
&\phi_{1,R}(u)=1-\rho(\frac{5}{7},\frac{6}{7})(\frac{u}{R}),&\quad &\psi_{2,R}=\rho(\frac{3}{7},\frac{4}{7})(\frac{u}{R}),\\
&\phi_{2,R}=\rho(\frac{1}{7},\frac{2}{7})(\frac{u}{R}),&\quad &\psi_{1,R}=1-\psi_{2,R}.
\end{aligned}\end{align}
We extend the functions in \eqref{e.97} symmetrically to the whole real line $(-\infty,\infty)$, then naturally extend to $Y\times(-\infty,\infty)$. Since these functions are constant outside $Y_{[-R,R]}$, we extend them trivially to the manifolds $M_R,\,M_{R,1},\,M_{R,2}$.

Let $\cE_{R,m}(t;x,y)$ denote the kernel of the operator $\exp(-tA_{R,m}^2)$. Let $\cE_{i,R}(t;x,y)$
be the kernel of the operator $\exp(-tA_{i,R,m}^2)$. Let $\cE_c(t;x,y)$ be the kernel of $e^{-t(D^2_{Y_{\R}}+m^2)}$ on the infinite cylinder $Y\times [-\infty,\infty]$. Let $\cE_{c,2}(t;x,y)$ be the kernel of $e^{-t(D_{Y_{\R^+}}-m\GS)_{V_+}^2}$ on $Y_{\R^+}$ with the boundary condition \eqref{e.126}. Let $\cE_{c,1}(t;(u,x),(v,y))$ be the kernel of $e^{-t(D_{Y_{\R^-}}-m\GS)_{V_-}^2}$ on $Y_{\R^-}$ with the boundary condition \eqref{e.195}.

We define a parametrix for $\cE_R(t;x,y)$ as follows: (cf. \cite{APS}, \cite[(4.2)]{DouWoj91})
\begin{align}\begin{aligned}\label{e.98}
Q_{R,m}(t;x,y):=\phi_{1,R}(x)\cE_c(t;x,y)\psi_{1,R}(y)+\phi_{2,R}(x)\cE_{R,m}(t;x,y)\psi_{2,R}(y).
\end{aligned}\end{align}
Similarly, we define parametrices for $\cE_{i,R}(t;x,y),\,i=1,2$,
\begin{align}\begin{aligned}\label{e.99}
Q_{i,R,m}(t;x,y):=\phi_{1,R}(x)\cE_{c,i}(t;x,y)\psi_{1,R}(y)+\phi_{2,R}(x)\cE_{R,m}(t;x,y)\psi_{2,R}(y),
\end{aligned}\end{align}
where $(x,y)\in M_{R,i}\times M_{R,i}$.

By Duhamel's principle, we have
\begin{align}\begin{aligned}\label{e.100}
\cE_{R,m}(t;x,y)=Q_{R,m}(t;x,y)+(\cE_{R,m}*\cC_{R,m})(t;x,y),
\end{aligned}\end{align}
where the error term $\cC_R(t;x,y)$ is defined as
\begin{align}\begin{aligned}\label{e.101}
\cC_{R,m}(t;x,y):=\left(\partial_t+A_{R,m}^2\right)Q_{R,m}(t;x,y),
\end{aligned}\end{align}
and
\begin{align}\begin{aligned}\label{e.102}
(\cE_{R,m}*\cC_{R,m})(t;x,y):=\int^t_0\int_{M_R}dz \cE_{R,m}(s;,x,z)\cC_{R,m}(t-s;z,y).
\end{aligned}\end{align}

Similarly, we define
\begin{align}\begin{aligned}\label{e.103}
\cC_{i,{R,m}}(t;x,y):=\left(\partial_t+A_{i,R,m}^2\right)Q_{i,R,m}(t;x,y),
\end{aligned}\end{align}
then we also have
\begin{align}\begin{aligned}\label{e.104}
\cE_{i,R,m}(t;x,y)=Q_{i,R,m}(t;x,y)+(\cE_{i,R,m}*\cC_{i,R,m})(t;x,y).
\end{aligned}\end{align}

\begin{prop}\label{p.3}We have
  \begin{align}\begin{aligned}\label{e.105}
\lim_{R\ra \infty}\lim_{s\ra 0}\frac{1}{\G((s+1)/2)}\int_0^{R^{2-\var}} t^{(s-1)/2}\int_{M_R}\tr\big[
(A_{R,m})_x(\cE_{R,m}*\cC_{R,m})(t;x,y)\big]_{x=y}dxdt=0.
\end{aligned}\end{align}
And similarly for $i=1,2$
\begin{align}\begin{aligned}\label{e.106}
\lim_{R\ra \infty}\lim_{s\ra 0}\frac{1}{\G((s+1)/2)}\int_0^{R^{2-\var}} t^{(s-1)/2}\int_{M_{R,i}}\tr\big[
(A_{i,R,m})_x(\cE_{i,R,m}*\cC_{i,R,m})(t;x,y)\big]_{x=y}dxdt=0.
\end{aligned}\end{align}
\end{prop}
\begin{proof}
Proposition \ref{p.3} follows by a similar argument as in \cite[(4.10)]{DouWoj91}, \cite[Lemma 3.2]{Wojcie} or \cite[Lemma 2.16]{Zhu16}.
\end{proof}

By Proposition \ref{p.3}, \eqref{e.100} and \eqref{e.104}, we see that the main contribution to $\delta^S_{m,R}(s)$ comes from $Q_{R}(t;x,y),\,Q_{i,R}(t;x,y)$ for $i=1,2$ in the adiabatic limit.
By our construction of the parametrices in \eqref{e.98}, \eqref{e.99}, the interior contributions cancel with each other, so we only need to consider the contributions from the cylinder parts. Next, we aim to compute the following term: for $x\in \cY{-R,R}$
\begin{align}\begin{aligned}\label{e.108}
\Delta_{c,R,m}(t;x):=\tr\Big[
(D-m\GS)_x \phi_{1,R}(x)(\cE_c-\sum_{i=1}^2\cE_{c,i})(t;x,y)\psi_{1,R}(y)\Big]_{x=y},
\end{aligned}\end{align}
where $D=\gamma(\partial_u+D_Y)$ and we have used the fact that both $A_{R,m}$ and $A_{i,R,m},\,i=1,2$ are of form $D-m\GS$ on the cylinder part $\cY{-R,R}$ (see \eqref{e.76}). Then by Proposition \ref{p.3}, \eqref{e.98}, \eqref{e.99}, \eqref{e.100}, \eqref{e.104} and \eqref{e.108}, we have
\begin{align}\begin{aligned}\label{e.109}
\lim_{R\ra \infty}\lim_{s\ra 0}\delta^S_{m,R}(s)=\lim_{R\ra \infty}\lim_{s\ra 0}\frac{1}{\G((s+1)/2)}\int_0^{R^{2-\var}} t^{(s-1)/2}\int_{\cY{-R,R}}\Delta_{c,R,m}(t;x)dxdt.
\end{aligned}\end{align}
By the definition of cut-off functions in \eqref{e.97} and \eqref{e.108}, in fact we have
\begin{align}\begin{aligned}\label{e.110}
\Delta_{c,R,m}(t;x):&=\tr\big[\phi_{1,R}(x)
(D-m\GS)_x (\cE_c-\sum_{i=1}^2\cE_{c,i})(t;x,y)\psi_{1,R}(y)\big]_{x=y}\\
&=\psi_{1,R}(x)\tr\big[
(D-m\GS)_x (\cE_c-\sum_{i=1}^2\cE_{c,i})(t;x,y)\big]_{x=y}.
\end{aligned}\end{align}

The heat kernels on the infinite or half-infinite cylinders has been explicitly calculated as in \eqref{e.142} and \eqref{e.160}.
Assume that $\set{\phi_{\lambda_i},\, i\in \N}$ is an orthonormal basis of ${\rm Ran}(P_>)$ consisting of the eigensections of $D_Y$ with eigenvalues $0<\lambda_1\leq \lambda_2\leq \cdots$ and $\set{\psi_i,\,1\leq i\leq n_+}$ is an orthonormal basis of $V_+$. For $(u,x),(v,y)\in Y_{\R^+}$, the kernel $\cE_{c,2}(t;(u,x),(v,y))$ of $e^{-t(D_{Y_{\R^+}}-m\GS)_{V_+}^2}$ with the boundary condition \eqref{e.126} is given by (see \eqref{e.142})
\begin{align}\begin{aligned}\label{e.143}
&\cE_{c,2}(t;(u,x),(v,y))\\
=&e^{-m^2t}\sum_{\lambda_k>0}\frac{e^{-\lambda_k^2t}}{\sqrt{4\pi t}}\left\{\exp\left(\frac{-(u-v)^2}{4t}\right)-\exp\left(\frac{-(u+v)^2}{4t}\right)\right\}
\phi_{\lambda_k}(x)\otimes \overline{\phi_{\lambda_k}}(y)\\
&+e^{-m^2t}\sum_{\lambda_k>0}\left\{\frac{e^{-\lambda_k^2t}}{\sqrt{4\pi t}}\left(\exp\left(\frac{-(u-v)^2}{4t}\right)+\exp\left(\frac{-(u+v)^2}{4t}\right)\right)\right.\\
&\left.\qquad\qquad\qquad\qquad-\lambda_k e^{\lambda_k(u+v)}\erfc\left(\frac{u+v}{2\sqrt{t}}+\lambda_k
\sqrt{t}\right)\right\}\gamma\phi_{\lambda_k}(x)\otimes \overline{\gamma\phi_{\lambda_k}}(y)
\\
&+\frac{e^{-m^2t}}{\sqrt{4\pi t}}\left\{\exp\left(\frac{-(u-v)^2}{4t}\right)-\exp\left(\frac{-(u+v)^2}{4t}\right)\right\}\sum_{i=1}^{n_+}
\psi_i(x)\otimes \overline{\psi_i}(y)\\
&+\frac{e^{-m^2t}}{\sqrt{4\pi t}}\left\{\exp\left(\frac{-(u-v)^2}{4t}\right)+\exp\left(\frac{-(u+v)^2}{4t}\right)\right\}\sum_{i=1}^{n_+}\gamma\psi_i(x)\otimes \overline{\gamma\psi_i}(y),
\end{aligned}\end{align}
where $\erfc(x)$ is the complementary error function defined in \eqref{e.144}.
Similarly, we can get the expression of the kernel $\cE_{c,1}(t;(u,x),(v,y))$ of $e^{-t(D_{Y_{\R^-}}-m\GS)_{V_-}^2}$ with the boundary condition \eqref{e.160},
\begin{align}\begin{aligned}\label{e.114}
&\cE_{c,1}(t;(u,x),(v,y))\\
=&e^{-m^2t}\sum_{\lambda_k>0}\frac{e^{-\lambda_k^2t}}{\sqrt{4\pi t}}\left\{\exp\left(\frac{-(u-v)^2}{4t}\right)-\exp\left(\frac{-(u+v)^2}{4t}\right)\right\}
\phi_{\lambda_k}(x)\otimes \overline{\phi_{\lambda_k}}(y)\\
&+e^{-m^2t}\sum_{\lambda_k>0}\left\{\frac{e^{-\lambda_k^2t}}{\sqrt{4\pi t}}\left(\exp\left(\frac{-(u-v)^2}{4t}\right)+\exp\left(\frac{-(u+v)^2}{4t}\right)\right)\right.\\
&\left.\qquad\qquad\qquad\qquad-\lambda_k e^{-\lambda_k(u+v)}\erfc\left(-\frac{u+v}{2\sqrt{t}}+\lambda_k
\sqrt{t}\right)\right\}\gamma\phi_{\lambda_k}(x)\otimes \overline{\gamma\phi_{\lambda_k}}(y)
\\
&+\frac{e^{-m^2t}}{\sqrt{4\pi t}}\left\{\exp\left(\frac{-(u-v)^2}{4t}\right)-\exp\left(\frac{-(u+v)^2}{4t}\right)\right\}\sum_{i=1}^{n_+}
\psi_i(x)\otimes \overline{\psi_i}(y)\\
&+\frac{e^{-m^2t}}{\sqrt{4\pi t}}\left\{\exp\left(\frac{-(u-v)^2}{4t}\right)+\exp\left(\frac{-(u+v)^2}{4t}\right)\right\}\sum_{i=1}^{n_+}\gamma\psi_i(x)\otimes \overline{\gamma\psi_i}(y).
\end{aligned}\end{align}

Moreover, the kernel $\cE_{c}(t;(u,x),(v,y))$ of $e^{-t(D^2_{Y_{\R}}+m^2)}$ is given by
\begin{align}\begin{aligned}\label{e.115}
&\cE_{c}(t;(u,x),(v,y))\\
=&\frac{e^{-m^2t}}{\sqrt{4\pi t}}\exp\Big(-\frac{(u-v)^2}{4t}\Big)\cdot e^{-tD^2_Y}(x,y)\\
=&\frac{e^{-m^2t}}{\sqrt{4\pi t}}\exp\Big(-\frac{(u-v)^2}{4t}\Big)\sum_{i=1}^\infty e^{-t\lambda^2_i}
\left(\gamma\phi_{\lambda_i}(x)\otimes \overline{\gamma\phi_{\lambda_i}(y)}+\phi_{\lambda_i}(x)\otimes \overline{\phi_{\lambda_i}(y)}\right)\\
&+\frac{e^{-m^2t}}{\sqrt{4\pi t}}\exp\Big(-\frac{(u-v)^2}{4t}\Big)\sum_{i=1}^{n_+}
\left(\gamma\psi_i(x)\otimes \overline{\gamma\psi_i(y)}+\psi_i(x)\otimes \overline{\psi_i(y)}\right).
\end{aligned}\end{align}

By \eqref{e.110}, we have for $x=(u,y)\in [-R,R]\times Y$
\begin{align}\begin{aligned}\label{e.116}
&\int_{\cY{-R,R}}\Delta_{c,R,m}(t;x)dx\\
=&\int_{-R}^R \int_Y \psi_{1,R}(u)\tr\big[
(D-m\GS)_x (\cE_c-\sum_{i=1}^2\cE_{c,i})(t;x,z)\big]_{x=z}dydu\\
=& \int_{-R}^R \psi_{1,R}(u)\int_Y  \tr\big[
D_x (\cE_c-\sum_{i=1}^2\cE_{c,i})(t;x,z)\big]_{x=z}dydu\\
&\qquad -m\int_{-R}^R\psi_{1,R}(u) \int_Y \tr\big[
(\GS)_x (\cE_c-\sum_{i=1}^2\cE_{c,i})(t;x,z)\big]_{x=z}dydu\\
=&-m\int_{-R}^R\psi_{1,R}(u) \int_Y \tr\big[
(\GS)_x (\cE_c-\sum_{i=1}^2\cE_{c,i})(t;x,z)\big]_{x=z}dydu,
\end{aligned}\end{align}
where for the last equality we have used the facts
\begin{align}\begin{aligned}\label{e.117}
\int_Y  \tr\big[
D_x \cE_c(t;x,y)\big]_{x=y}dy=0,\quad \int_Y  \tr\big[
D_x\cE_{c,i}(t;x,y)\big]_{x=y}dy=0, \quad \text{for}\,\,i=1,2,
\end{aligned}\end{align}
since we have $
\iprod{\psi_i,\gamma\psi_i}_Y=0,\, \iprod{\phi_{\lambda_k},\gamma\phi_{\lambda_k}}_Y=0$.
By \eqref{e.39} and \eqref{e.114}, we have
\begin{align}\begin{aligned}\label{e.146}
&\int_Y \tr\big[
(\GS)_x\cE_{c,1})(t;x,z)\big]_{x=z}dy\\
=&e^{-m^2t}\sum_{\lambda_k>0}\frac{e^{-\lambda_k^2t}}{\sqrt{4\pi t}}\left(1-e^{-u^2/t}\right)
\iprod{\GS\phi_{\lambda_k}, \phi_{\lambda_k}}_Y\\
&+e^{-m^2t}\sum_{\lambda_k>0}\left\{\frac{e^{-\lambda_k^2t}}{\sqrt{4\pi t}}\left(1+e^{-u^2/t}\right)
-\lambda_k e^{-2\lambda_ku}\erfc\left(-\frac{u}{\sqrt{t}}+\lambda_k
\sqrt{t}\right)\right\}\iprod{\GS\gamma\phi_{\lambda_k}, \gamma\phi_{\lambda_k}}_Y
\\
&+\frac{e^{-m^2t}}{\sqrt{4\pi t}}\left\{1-e^{-u^2/t}\right\}\sum_{i=1}^{n_+}
\iprod{\GS\psi_i, \psi_i}_Y+\frac{e^{-m^2t}}{\sqrt{4\pi t}}\left\{1+e^{-u^2/t}\right\}\sum_{i=1}^{n_+}\iprod{\GS\gamma\psi_i, \gamma\psi_i}_Y\\
=&e^{-m^2t}\sum_{\lambda_k>0}\frac{e^{-\lambda_k^2t}}{\sqrt{4\pi t}}\left(1-e^{-u^2/t}\right)
\iprod{\GS\phi_{\lambda_k},\phi_{\lambda_k}}_Y\\
&+e^{-m^2t}\sum_{\lambda_k>0}\left\{\frac{e^{-\lambda_k^2t}}{\sqrt{4\pi t}}\left(-1-e^{-u^2/t}\right)
+\lambda_k e^{-2\lambda_ku}\erfc\left(-\frac{u}{\sqrt{t}}+\lambda_k
\sqrt{t}\right)\right\}\iprod{\GS\phi_{\lambda_k}, \phi_{\lambda_k}}_Y\\
&+\frac{e^{-m^2t}}{\sqrt{4\pi t}}\left\{1-e^{-u^2/t}\right\}\sum_{i=1}^{n_+}
\norm{\psi_i}^2_Y-\frac{e^{-m^2t}}{\sqrt{4\pi t}}\left\{1+e^{-u^2/t}\right\}\sum_{i=1}^{n_+}\norm{\psi_i}^2_Y\\
=&e^{-m^2t}\sum_{\lambda_k>0}\left\{-2\frac{e^{-\lambda_k^2t-u^2/t}}{\sqrt{4\pi t}}
+\lambda_k e^{-2\lambda_ku}\erfc\left(-\frac{u}{\sqrt{t}}+\lambda_k
\sqrt{t}\right)\right\}\iprod{\GS\phi_{\lambda_k}, \phi_{\lambda_k}}_Y\\
&\qquad\qquad\qquad\qquad\qquad-2\frac{e^{-m^2t-u^2/t}}{\sqrt{4\pi t}}\sum_{i=1}^{n_+}
\norm{\psi_i}^2_Y,
\end{aligned}\end{align}
where we have used
\begin{align}\begin{aligned}\label{e.274}
\iprod{\GS\gamma\psi,\gamma\psi}_Y&=\iprod{-\gamma\GS\psi,\gamma\psi}_Y=-\iprod{\GS\psi,\gamma^*\gamma\psi}_Y\\
&=-\iprod{\GS\psi,-\gamma^2\psi}_Y=-\iprod{\GS\psi,\psi}_Y\\
&=-\iprod{\psi,\psi}_Y=-\norm{\psi}_Y^2.
\end{aligned}\end{align}

By \eqref{e.39}, \eqref{e.146} and $\norm{\psi_i}_Y=1,\,1\leq i\leq n_+$, we have
\begin{align}\begin{aligned}\label{e.118}
&\int_{-R}^0\psi_{1,R}(u) \int_Y \tr\big[
(\GS)_x\cE_{c,1})(t;x,z)\big]_{x=z}dydu\\
=&e^{-m^2t}\int_{-R}^0\psi_{1,R}(u)\sum_{\lambda_k>0}\left\{-\frac{2e^{-\lambda_k^2t-u^2/t}}{\sqrt{4\pi t}}
+\lambda_k e^{-2\lambda_ku}\erfc\left(-\frac{u}{\sqrt{t}}+\lambda_k
\sqrt{t}\right)\right\}du\iprod{\GS\phi_{\lambda_k}, \phi_{\lambda_k}}_Y\\
&\qquad\qquad\qquad\qquad\qquad-2n_+\int_{-R}^0\psi_{1,R}(u)\frac{e^{-m^2t-u^2/t}}{\sqrt{4\pi t}}du.
\end{aligned}\end{align}

By \eqref{e.39} and \eqref{e.143}, we get
\begin{align}\begin{aligned}\label{e.148}
&\int_Y \tr\big[
(\GS)_x\cE_{c,2})(t;x,z)\big]_{x=z}dy\\
=&e^{-m^2t}\sum_{\lambda_k>0}\frac{e^{-\lambda_k^2t}}{\sqrt{4\pi t}}\left(1-e^{-u^2/t}\right)
\iprod{\GS\phi_{\lambda_k},\phi_{\lambda_k}}_Y\\
&+e^{-m^2t}\sum_{\lambda_k>0}\left\{\frac{e^{-\lambda_k^2t}}{\sqrt{4\pi t}}\left(1+e^{-u^2/t}\right)-\lambda_k e^{2\lambda_ku}\erfc\left(\frac{u}{\sqrt{t}}+\lambda_k
\sqrt{t}\right)\right\}\iprod{\GS\gamma\phi_{\lambda_k}, \gamma\phi_{\lambda_k}}_Y
\\
&+\frac{e^{-m^2t}}{\sqrt{4\pi t}}\left(1-e^{-u^2/t}\right)
\sum_{i=1}^{n_+}\iprod{\GS\psi_i,\psi_i}_Y+\frac{e^{-m^2t}}{\sqrt{4\pi t}}\left(1+e^{-u^2/t}\right)\sum_{i=1}^{n_+}\iprod{\GS\gamma\psi_i, \gamma\psi_i}_Y\\
=
&e^{-m^2t}\sum_{\lambda_k>0}\left\{-2\frac{e^{-\lambda_k^2t-u^2/t}}{\sqrt{4\pi t}}+\lambda_k e^{2\lambda_ku}\erfc\left(\frac{u}{\sqrt{t}}+\lambda_k
\sqrt{t}\right)\right\}\iprod{\GS\phi_{\lambda_k}, \phi_{\lambda_k}}_Y\\
&\qquad\qquad\qquad\qquad\qquad-2\frac{e^{-m^2t-u^2/t}}{\sqrt{4\pi t}}
\sum_{i=1}^{n_+}\norm{\psi_i}^2_Y.
\end{aligned}\end{align}

We get by \eqref{e.148} and $\norm{\psi_i}_Y=1,\,1\leq i\leq n_+$,
\begin{align}\begin{aligned}\label{e.119}
&\int_{0}^R\psi_{1,R}(u) \int_Y \tr\big[
(\GS)_x\cE_{c,2})(t;x,z)\big]_{x=z}dydu\\
=&e^{-m^2t}\int_{0}^R\psi_{1,R}(u)\sum_{\lambda_k>0}\left\{\frac{-2e^{-u^2/t-\lambda_k^2t}}{\sqrt{4\pi t}}+\lambda_k e^{2\lambda_ku}\erfc\left(\frac{u}{\sqrt{t}}+\lambda_k
\sqrt{t}\right)\right\}du\iprod{\GS\phi_{\lambda_k}, \phi_{\lambda_k}}_Y\\
&\qquad\qquad\qquad\qquad\qquad-2n_+\int_{0}^R\psi_{1,R}(u)\frac{e^{-m^2t-u^2/t}}{\sqrt{4\pi t}}
du.
\end{aligned}\end{align}
By \eqref{e.39} and \eqref{e.115}, we have
\begin{align}\begin{aligned}\label{e.150}
&\int_Y \tr\big[
(\GS)_x\cE_{c})(t;x,z)\big]_{x=z}dy\\
=&\frac{e^{-m^2t}}{\sqrt{4\pi t}}\sum_{i=1}^\infty e^{-t\lambda^2_i}
\left(\iprod{\GS\gamma\phi_{\lambda_i}, \gamma\phi_{\lambda_i}}+\iprod{\GS\phi_{\lambda_i}, \phi_{\lambda_i}}\right)\\
&+\frac{e^{-m^2t}}{\sqrt{4\pi t}}\sum_{i=1}^{n_+}
\left(\iprod{\GS\gamma\psi_i, \gamma\psi_i}+\iprod{\GS\psi_i, \psi_i }\right)\\
=&\frac{e^{-m^2t}}{\sqrt{4\pi t}}\sum_{i=1}^\infty e^{-t\lambda^2_i}
\left(-\iprod{\GS\phi_{\lambda_i}, \phi_{\lambda_i}}+\iprod{\GS\phi_{\lambda_i}, \phi_{\lambda_i}}\right)\\
&+\frac{e^{-m^2t}}{\sqrt{4\pi t}}\sum_{i=1}^{n_+}
\left(-\iprod{\GS\psi_i, \psi_i}+\iprod{\GS\psi_i, \psi_i }\right)=0.
\end{aligned}\end{align}
Hence we get
\begin{align}\begin{aligned}\label{e.120}
&\int_{-R}^R\psi_{1,R}(u) \int_Y \tr\big[
(\GS)_x\cE_{c})(t;x,z)\big]_{x=z}dydu=0.
\end{aligned}\end{align}

By \eqref{e.118}, \eqref{e.119} and \eqref{e.120}, we have
\begin{align}\begin{aligned}\label{e.151}
&-\int_{-R}^R\psi_{1,R}(u) \int_Y \tr\big[
(\GS)_x (\cE_c-\sum_{i=1}^2\cE_{c,i})(t;x,z)\big]_{x=z}dydu\\
=
&e^{-m^2t}\int_{-R}^0\psi_{1,R}(u)\sum_{\lambda_k>0}\left\{-\frac{2e^{-\lambda_k^2t-u^2/t}}{\sqrt{4\pi t}}
+\lambda_k e^{-2\lambda_ku}\erfc\left(-\frac{u}{\sqrt{t}}+\lambda_k
\sqrt{t}\right)\right\}du\iprod{\GS\phi_{\lambda_k}, \phi_{\lambda_k}}_Y
\\
+&e^{-m^2t}\int_{0}^R\psi_{1,R}(u)\sum_{\lambda_k>0}\left\{\frac{-2e^{-u^2/t-\lambda_k^2t}}{\sqrt{4\pi t}}+\lambda_k e^{2\lambda_ku}\erfc\left(\frac{u}{\sqrt{t}}+\lambda_k
\sqrt{t}\right)\right\}du\iprod{\GS\phi_{\lambda_k}, \phi_{\lambda_k}}_Y\\
&-2n_+\int_{-R}^0\psi_{1,R}(u)\frac{e^{-m^2t-u^2/t}}{\sqrt{4\pi t}}du -2n_+\int_{0}^R\psi_{1,R}(u)\frac{e^{-m^2t-u^2/t}}{\sqrt{4\pi t}}
du\\
=
&2e^{-m^2t}\int^{R}_0\psi_{1,R}(u)\sum_{\lambda_k>0}\left\{-\frac{2e^{-\lambda_k^2t-u^2/t}}{\sqrt{4\pi t}}
+\lambda_k e^{2\lambda_ku}\erfc\left(\frac{u}{\sqrt{t}}+\lambda_k
\sqrt{t}\right)\right\}du\iprod{\GS\phi_{\lambda_k}, \phi_{\lambda_k}}_Y\\
&\qquad\qquad\qquad\qquad\qquad\qquad -4n_+\int_{0}^R\psi_{1,R}(u)\frac{e^{-m^2t-u^2/t}}{\sqrt{4\pi t}}
du.
\end{aligned}\end{align}

By \eqref{e.116} and \eqref{e.151}, we get
\begin{align}\begin{aligned}\label{e.152}
&\int_{\cY{-R,R}}\Delta_{c,R,m}(t;x)dx\\
=&-m\int_{-R}^R\psi_{1,R}(u) \int_Y \tr\big[
(\GS)_x (\cE_c-\sum_{i=1}^2\cE_{c,i})(t;x,z)\big]_{x=z}dydu\\
=
&me^{-m^2t}\int^{R}_0\psi_{1,R}(u)\sum_{\lambda_k>0}\left\{-\frac{2e^{-\lambda_k^2t-u^2/t}}{\sqrt{\pi t}}
+2\lambda_k e^{2\lambda_ku}\erfc\left(\frac{u}{\sqrt{t}}+\lambda_k
\sqrt{t}\right)\right\}du\iprod{\GS\phi_{\lambda_k}, \phi_{\lambda_k}}_Y\\
&\qquad\qquad\qquad\qquad\qquad\qquad -4mn_+\int_{0}^R\psi_{1,R}(u)\frac{e^{-m^2t-u^2/t}}{\sqrt{4\pi t}}
du.
\end{aligned}\end{align}

Set
\begin{align}\begin{aligned}\label{e.169}
&a(R)=\int_0^R\psi_{1,R}(u)\frac{2e^{-u^2/t}}{\sqrt{\pi t}}du.
\end{aligned}\end{align}
Then we have
\begin{align}\begin{aligned}\label{e.170}
\lim_{R\ra \infty}a(R)=1.
\end{aligned}\end{align}
By \eqref{e.169}, we get
\begin{align}\begin{aligned}\label{e.171}
&\int_{\cY{-R,R}}\Delta_{c,R,m}(t;x)dx\\
=
&me^{-m^2t}\sum_{\lambda_k>0}\left\{\int^{R}_0\psi_{1,R}(u)2\lambda_k e^{2\lambda_ku} \erfc\left(\frac{u}{\sqrt{t}}+\lambda_k
\sqrt{t}\right)du-e^{-\lambda_k^2t}a(R)\right\}\iprod{\GS\phi_{\lambda_k}, \phi_{\lambda_k}}_Y\\
&\qquad\qquad\qquad\qquad\qquad\qquad -mn_+e^{-m^2t}a(R).
\end{aligned}\end{align}
Using integration by parts and \eqref{e.170}, we get
\begin{align}\begin{aligned}\label{e.172}
&\int^{R}_0\psi_{1,R}(u)2x e^{2xu} \erfc\left(\frac{u}{\sqrt{t}}+x
\sqrt{t}\right)du\\
=&- \erfc(x\sqrt{t})+e^{-x^2t}a(R)-\int^R_0\psi'_{1,R}(u)e^{2x u}\erfc\left(\frac{u}{\sqrt{t}}+x
\sqrt{t}\right)du.
\end{aligned}\end{align}
By \eqref{e.171} and \eqref{e.172}, we have
\begin{align}\begin{aligned}\label{e.173}
&\int_{\cY{-R,R}}\Delta_{c,R,m}(t;x)dx\\
=
&me^{-m^2t}\sum_{\lambda_k>0}\left\{ -\erfc(\lambda_k\sqrt{t})-\int^R_0\psi'_{1,R}(u)e^{2\lambda_k u}\erfc\left(\frac{u}{\sqrt{t}}+\lambda_k
\sqrt{t}\right)du\right\}\iprod{\GS\phi_{\lambda_k}, \phi_{\lambda_k}}_Y\\
&\qquad\qquad\qquad\qquad\qquad\qquad -mn_+e^{-m^2t}a(R).
\end{aligned}\end{align}

Set
\begin{align}\begin{aligned}\label{e.167}
G(R,t)&=\int^{R}_0\psi'_{1,R}(u)\sum_{\lambda_k>0}e^{2\lambda_ku} \erfc\left(\frac{u}{\sqrt{t}}+\lambda_k
\sqrt{t}\right)
du\iprod{\GS\phi_{\lambda_k}, \phi_{\lambda_k}}_Y,\\
\Theta(t)&=-me^{-m^2t}\sum_{\lambda_k>0}\erfc(\lambda_k\sqrt{t})\iprod{\GS\phi_{\lambda_k}, \phi_{\lambda_k}}_Y.
\end{aligned}\end{align}
By \eqref{e.173} and \eqref{e.167}, we have
\begin{align}\begin{aligned}\label{e.168}
\int_{\cY{-R,R}}\Delta_{c,R,m}(t;x)dx
=\Theta(t)-me^{-m^2t}\left(G(R,t)+n_+a(R)\right).
\end{aligned}\end{align}

Set
\begin{align}\begin{aligned}\label{e.178}
\xi_R(s;D_Y)=\frac{1}{\G((s+1)/2)}\int_0^{R^{2-\var}} t^{(s-1)/2}\Theta(t)dt.
\end{aligned}\end{align}

\begin{lemma}\label{l.1}
We have
  \begin{align}\begin{aligned}\label{e.174}
&\lim_{R\ra \infty}\lim_{s\ra 0}\frac{1}{\G((s+1)/2)}\int_0^{R^{2-\var}}t^{(s-1)/2}me^{-m^2t}G(R,t)dt=0,\\
&\lim_{R\ra \infty}\lim_{s\ra 0}\frac{1}{\G((s+1)/2)}\int_0^{R^{2-\var}}t^{(s-1)/2}me^{-m^2t}n_+a(R)dt=n_+.
\end{aligned}\end{align}
\end{lemma}
\begin{proof}
 Set $\psi_1(u)=\psi_{1,R}(u)|_{R=1}$. By \eqref{e.96} and \eqref{e.97}, we have
 \begin{align}\begin{aligned}\label{e.200}
{\rm supp }\left(\psi'_1\right)\subset [\frac{3}{7},\frac{4}{7}]\times Y.
\end{aligned}\end{align}
  Since $\erfc(r)\leq c e^{-r^2}$, $\psi_{1,R}(u)=\psi_1(\frac{u}{R})$ and $\norm{\phi_{\lambda_k}}^2_Y=1$, we get by \eqref{e.167} and \eqref{e.200}
\begin{align}\begin{aligned}\label{e.175}
\abs{G(R,t)}\leq& c\int^{\infty}_0\psi'_{1}(u/R)\frac{1}{R}\sum_{\lambda_k>0}e^{2\lambda_ku} \exp\left(-u^2/t-2\lambda_ku-\lambda^2_k
t\right)
du\norm{\phi_{\lambda_k}}^2_Y\\
\leq& c\int^{\infty}_0\psi'_{1}(u)\sum_{\lambda_k>0} \exp\left(-\frac{R^2u^2}{t}-\lambda^2_k
t\right)
du\\
\leq &c\int^{\frac{4}{7}}_{\frac{3}{7}}\psi'_{1}(u)e^{-\frac{R^2u^2}{t}}du\sum_{\lambda_k>0} e^{-\lambda^2_k
t}\leq ce^{-c'\frac{R^2}{t}}\sum_{\lambda_k>0} e^{-\lambda^2_k
t}.
\end{aligned}\end{align}
By \eqref{e.175}, we get
\begin{align}\begin{aligned}\label{e.177}
&\left|\frac{1}{\G((s+1)/2)}\int_0^{R^{2-\var}}t^{(s-1)/2}me^{-m^2t}G(R,t)dt\right|\\
\leq &c\frac{1}{\left|\G((s+1)/2)\right|}\int_0^{R^{2-\var}}t^{(s-1)/2}e^{-c'\frac{R^2}{t}}\sum_{\lambda_k>0} e^{-(\lambda_k^2+m^2)
t}dt\\
\leq&c\frac{1}{\left|\G((s+1)/2)\right|}\int_0^{\infty}t^{(s-1)/2}e^{-c'\frac{R^2}{t}}\tr{e^{-t(D^2_Y+m^2)}}dt.
\end{aligned}\end{align}
Then the first equation in \eqref{e.174} follows from \eqref{e.177}. By \eqref{e.170}, we get
\begin{align}\begin{aligned}\label{e.275}
&\lim_{R\ra \infty}\lim_{s\ra 0}\frac{1}{\G((s+1)/2)}\int_0^{R^{2-\var}}t^{(s-1)/2}me^{-m^2t}n_+a(R)dt\\
=&n_+\lim_{s\ra 0}\frac{m}{\G((s+1)/2)}\int_0^{\infty}t^{(s-1)/2}e^{-m^2t}dt\\
=&n_+\lim_{s\ra 0}m^{-s}=n_+.
\end{aligned}\end{align}
We get the second equation in \eqref{e.174} by \eqref{e.275}.
The proof of Lemma \ref{l.1} is completed.
\end{proof}

By Lemma \ref{l.1}, \eqref{e.168} and \eqref{e.178}, we get
\begin{align}\begin{aligned}\label{e.161}
\lim_{R\ra \infty}\lim_{s\ra 0}\delta^S_{m,R}(s)&=\lim_{R\ra \infty}\lim_{s\ra 0}\frac{1}{\G((s+1)/2)}\int_0^{R^{2-\var}} t^{(s-1)/2}\int_{\cY{-R,R}}\Delta_{c,R,m}(t;x)dxdt\\
&=\lim_{R\ra\infty}\lim_{s\ra 0}\xi_R(s;D_Y)-n_+.
\end{aligned}\end{align}

\begin{thm}\label{t.6}
For the small time contribution (see \eqref{e.85}), we have
\begin{align}\begin{aligned}\label{e.201}
\lim_{R\ra \infty}\lim_{s\ra 0}\delta^S_{m,R}(s)=\lim_{R\ra\infty}\lim_{s\ra 0}\xi_R(s;D_Y)-n_+.
\end{aligned}\end{align}
\end{thm}

\subsection{Asymptotic gluing formula}

By Proposition \ref{p.1}, Theorem \ref{t.4}, Theorem \ref{t.6}, \eqref{e.78}, \eqref{e.81}, \eqref{e.84}-\eqref{e.86}, we get the following asymptotic gluing formula for eta invariants of massive Dirac type operators.

\begin{thm}\label{t.7}
Under the assumption \ref{a.1}, for $m>0$ sufficiently large we have
  \begin{align}\begin{aligned}\label{e.202}
&\lim_{s\ra 0}\left(\eta(s,D_{M}+mf)-\eta(s,\Domm)-\eta(s,\Dtmp)\right)\\
=&\lim_{R\ra\infty}\lim_{s\ra 0}\xi_R(s;D_Y)-n_+.
\end{aligned}\end{align}
\end{thm}

\begin{rem}\label{r.1}
  Note that the right hand side of \eqref{e.202} depends only on the boundary Dirac operator $D_Y$ in our geometric setting (see \eqref{e.75}, \eqref{e.76},\eqref{e.77}). Hence to compute the limit $\lim_{R\ra\infty}\lim_{s\ra 0}\xi_R(s;D_Y)$, we only need to know the geometric information near the cutting hypersurface. The author believes that this limit exists, but does not know how to prove it.
However, we don't need to know whether this limit exists for our applications of the asymptotic gluing formula \eqref{e.202}. If the right hand side of Equation \eqref{e.202} is divergent, all the divergences that occur will eventually cancel out in our final formulas.
\end{rem}

\section{Proof of Theorems \ref{t.3}, \ref{t.5}}\label{s4}

In this section, we will apply the asymptotic gluing formula in Theorem \ref{t.7} to prove our main results Theorems \ref{t.3}, \ref{t.5}.

\subsection{Decomposition of the domain-wall eta invariants}
Recall that $M_+$ and $M_-$ are two copies of $M$ with cylinder ends $Y_{[0,1]}$ and $Y_{[-1,0]}$ respectively (see Figure \ref{figure1}). The double manifold is defined as
\begin{align}\begin{aligned}\label{e.203}
\bM:=M_-\cup_Y M_+,
\end{aligned}\end{align}
which has $Y_{[-1,1]}$ as its cylinder part.

\begin{figure}
 \centering
  \includegraphics[width=10cm]{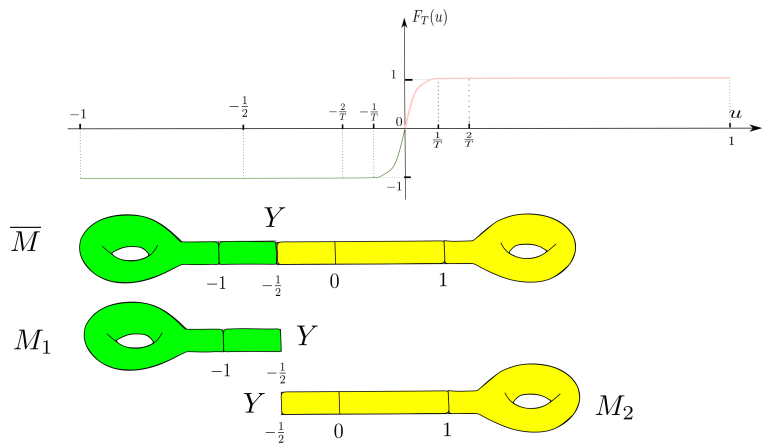}\\
\caption{Gluing setting for domain-wall massive Dirac operators. }\label{figure4}
\end{figure}

Now we cut $\bM$ into two pieces $M_1$,\, $M_2$ by using the hypersurface $Y_{\set{-\frac{1}{2}}}$ (see Figure \ref{figure4}), i.e.,
 \begin{align}\begin{aligned}\label{e.204}
\bM:=M_1\cup_{Y_{\set{-\onehalf}}} M_2,
\end{aligned}\end{align}
 such $M_1$ has the collar neighborhood $Y_{[-1,-\onehalf]}$ and $M_2$ has the collar neighborhood $Y_{[-\onehalf,1]}$.

Given $T>0$, let $F_T(u):[-R,1]\ra [-1,1]$ be a smooth increasing function, such that (see Figure \ref{figure4})
\begin{align}\begin{aligned}\label{e.52}
F_T(u)=-1,\quad \text{for}\, u\leq -\frac{1}{T}\quad ;F_T(u)=1,\quad \text{for}\,u\geq \frac{1}{T}.
\end{aligned}\end{align}
The function $F_T(u)$ can be naturally extended to the manifolds $\bM,\,M_{1},\,M_{2}$.
Note that $F_T(u)$ are constant functions, with value equal to $-1$, on the manifolds $\,M_{1}$.

In this case, the smooth section $f$ of $\End(S)$ used in Assumption \ref{a.1} is equal to $F_T(u)\GS$. For $m>0$ sufficiently large, we have
\begin{align}\begin{aligned}\label{e.205}
&\Ker_{L^2}(D_{M_{\infty,1}}+mF_T(u)\GS)=0,\quad \Ker_{L^2}(D_{M_{\infty,2}}+mF_T(u)\GS)=0,\\
&\quad \Ker(D_Y+m\gamma\GS)=\Ker(D^2_Y+m^2)=0,
\end{aligned}\end{align}
where we extend $F_T(u)\GS$ to the manifolds $M_{\infty,1},\,M_{\infty,2}$ with infinity cylinder ends in an obvious way. Hence Assumption \ref{a.1} is satisfied. Equation \eqref{e.76} is also satisfied by our definition of $F_T(u)$ in \eqref{e.52}.

Recall that $V_{\pm}\subset V=\Ker D_Y$ be the orthogonal Lagrangian subspaces induced in \eqref{e.70}.
We impose the generalized Atiyah-Patodi-Singer boundary condition with respect to the spectral projection $\Pi_{V_-}$  (see \eqref{e.124}) on $M_1$
and the one with respect to $\Pi_{V_+}$ on $M_2$.

By the asymptotic gluing formula of eta invariants in Theorem \ref{t.7}, we have for $m>0$ sufficiently large
\begin{align}\begin{aligned}\label{e.206}
&\lim_{s\ra 0}\left(\eta(s,D_{\bM}+mF_T(u))-\eta(s,(D_{M_1}-m\GS)_{V_-})-\eta(s,(D_{M_2}+mF_T(u))_{V_+})\right)\\
=&\lim_{R\ra\infty}\lim_{s\ra 0}\xi_R(s;D_Y)-n_+.
\end{aligned}\end{align}

Similarly applying Theorem \ref{t.7} to Pauli-Villars massive Dirac operators on $\bM$ (see \eqref{e.23}), we have for $m>0$ sufficiently large
\begin{align}\begin{aligned}\label{e.207}
&\lim_{s\ra 0}\left(\eta(s,D_{\bM}-m\GS)-\eta(s,(D_{M_1}-m\GS)_{V_-})-\eta(s,(D_{M_2}-m\GS)_{V_+})\right)\\
=&\lim_{R\ra\infty}\lim_{s\ra 0}\xi_R(s;D_Y)-n_+.
\end{aligned}\end{align}
By Remark \ref{r.1}, the right hand sides of \eqref{e.206} and \eqref{e.207} are the same, and will cancel out. Then by \eqref{e.206}, \eqref{e.207}, we obtain the following theorem.

\begin{thm}\label{t.8}
  For $m>0$ sufficiently large, we have the following identity of eta invariants,
 \begin{align}\begin{aligned}\label{e.208}
&\eta(D_{\bM}+mF_T\GS)-\eta(D_{\bM}-m\GS)\\
=&\eta((D_{M_2}+mF_T\GS)_{V_+})-\eta((D_{M_2}-m\GS)_{V_+}).
\end{aligned}\end{align}
\end{thm}

\begin{figure}
 \centering
  \includegraphics[width=8cm]{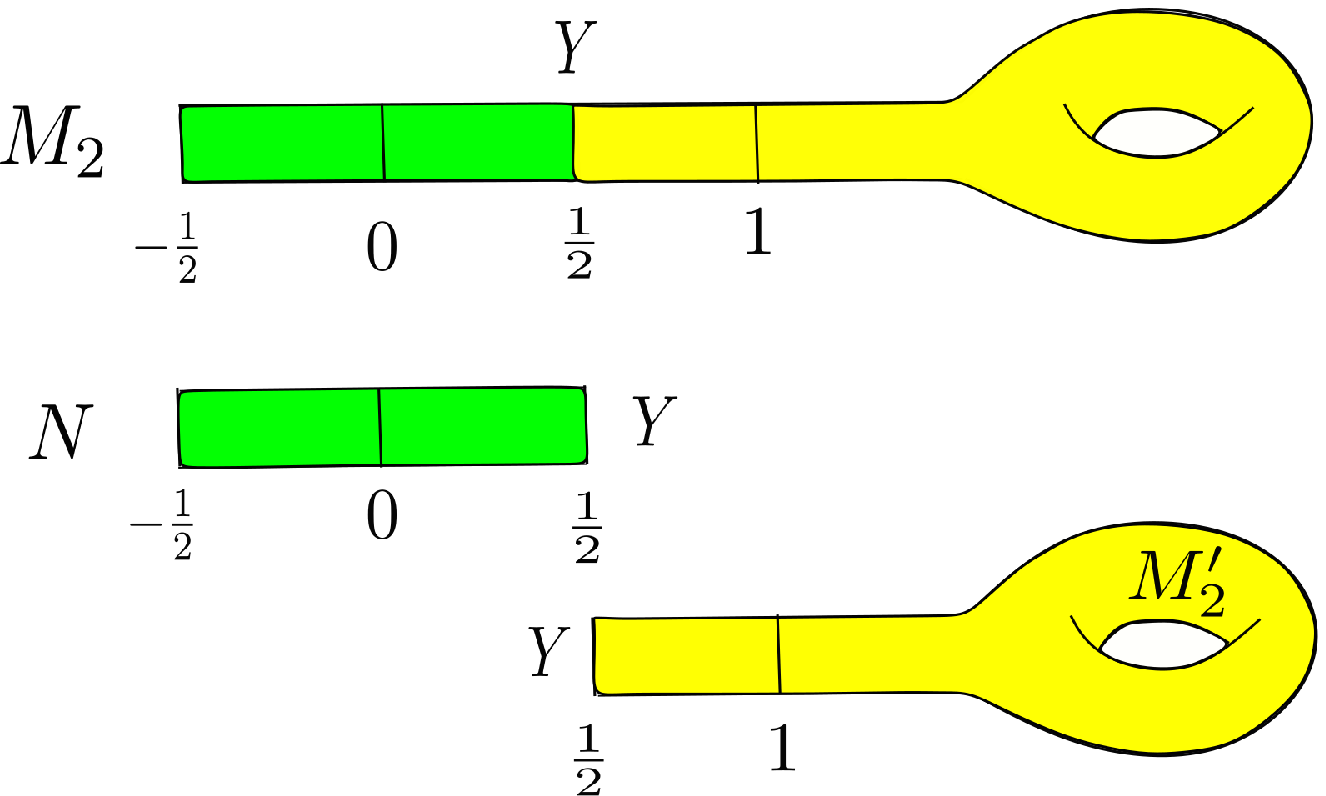}\\
\caption{Gluing setting for $M_2=N\cup_{Y_{\set{\onehalf}}}M'_2$. }\label{figure6}
\end{figure}

Note that $M_2$ has the cylinder end $Y\times [-\onehalf,1]$. We cut $M_2$ into two pieces $N$ and $M'_2$, i.e.,
\begin{align}\begin{aligned}\label{e.215}
M_2=N\cup_{Y_{\set{\onehalf}}}M'_2,
\end{aligned}\end{align}
such that $N=Y\times [-\onehalf,\onehalf]$ and $M'_2\subset M_2$ has the cylinder end $Y\times [\onehalf,1]$ (see Figure \ref{figure6}).

By using the asymptotic gluing formula \eqref{e.202} in Theorem \ref{t.7} and the same argument as above, we get by \eqref{e.52} and \eqref{e.215}
\begin{align}\begin{aligned}\label{e.216}
&\eta((D_{M_2}+mF_T\GS)_{V_+})-\eta((D_{M_2}+m\GS)_{V_+})\\
=&\eta((D_{M'_2}+mF_T\GS)_{V_+})-\eta((D_{M'_2}+m\GS)_{V_+})\\
&\quad\quad +\eta((D_{N}+mF_T\GS)_{V_+,V_-})-\eta((D_{N}+m\GS)_{V_+,V_-})\\
=&\eta((D_{M'_2}+m\GS)_{V_+})-\eta((D_{M'_2}+m\GS)_{V_+})\\
&\quad\quad +\eta((D_{N}+mF_T\GS)_{V_+,V_-})-\eta((D_{N}+m\GS)_{V_+,V_-})\\
=&\eta((D_{N}+mF_T\GS)_{V_+,V_-})-\eta((D_{N}+m\GS)_{V_+,V_-}),
\end{aligned}\end{align}
where $\eta((D_{N}+mF_T\GS)_{V_+,V_-}),\,\eta((D_{N}+m\GS)_{V_+,V_-})$ denote the eta invariants associated to massive Dirac operators over $N\cong Y\times [-\onehalf,\onehalf]$ with the generalized Atiyah-Patodi-Singer boundary conditions given by $\Pi_{V_+}$ at $Y_{\set{-\onehalf}}$ and by $\Pi_{V_-}$ at $Y_{\set{\onehalf}}$.

By \eqref{e.208}, \eqref{e.216}, we have for $m>0$ sufficiently large
\begin{align}\begin{aligned}\label{e.217}
&\eta(D_{\bM}+mF_T\GS)-\eta(D_{\bM}-m\GS)\\
=&\eta((D_{M_2}+mF_T\GS)_{V_+})-\eta((D_{M_2}+m\GS)_{V_+})\\
&\qquad\qquad +\eta((D_{M_2}+m\GS)_{V_+})-\eta((D_{M_2}-m\GS)_{V_+})\\
=&\eta((D_{M_2}+m\GS)_{V_+})-\eta((D_{M_2}-m\GS)_{V_+})\\
&\quad\quad +\eta((D_{N}+mF_T\GS)_{V_+,V_-})-\eta((D_{N}+m\GS)_{V_+,V_-}).
\end{aligned}\end{align}

By Definition \ref{d.2} and \eqref{e.32}, the right hand side of \eqref{e.47} in Theorem \ref{t.3} is just one half of the left hand side of \eqref{e.217}. To prove Theorem \ref{t.3}, we will show that the first two terms at the right hand side of \eqref{e.217} gives the generalized Atiyah-Patodi-Singer index. In next subsection, we will handle the last two terms, i.e.,
\begin{align}\begin{aligned}\label{e.229}
\eta((D_{N}+mF_T\GS)_{V_+,V_-})-\eta((D_{N}+m\GS)_{V_+,V_-})
\end{aligned}\end{align}
at the right hand side of \eqref{e.217} which are contributions from the finite cylinder $N=Y\times [-\onehalf,\onehalf]$.

\subsection{Eta invariants on the finite cylinder}
In this subsection, we will compute the difference of the eta invariants in \eqref{e.229} on the finite cylinder $N\cong Y\times[-\onehalf,\onehalf]$.

\begin{figure}
 \centering
  \includegraphics[width=10cm]{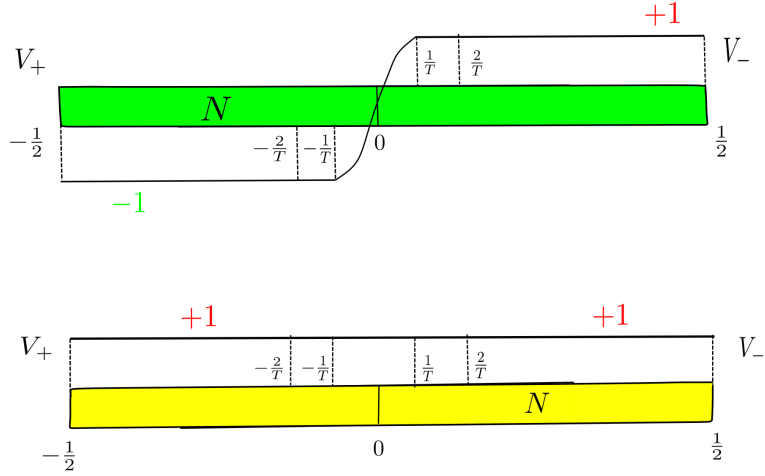}\\
\caption{Eta invariants on finite cylinders }\label{figure3}
\end{figure}

\begin{figure}
 \centering
  \includegraphics[width=8cm]{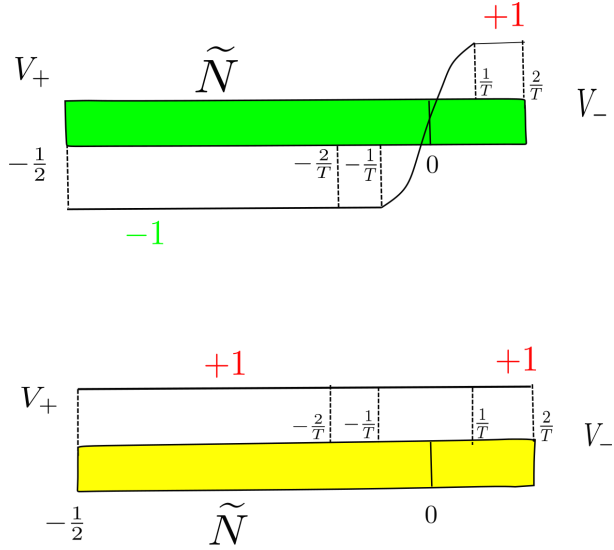}\\
\caption{Deform the cylinders }\label{figure9}
\end{figure}

We handle the terms in \eqref{e.229} by using a deformation argument.
 Recall that $\eta((D_N+mF_T\GS)_{V_+,V_-})$, $\eta((D_{N}+m\GS)_{V_+,V_-})$  are the eta invariants respectively associated with the Dirac operators $(D_N+mF_T\GS)_{V_+,V_-}$, $(D_N+m\GS)_{V_+,V_-}$ with the following boundary conditions (see Figure \ref{figure3}),
\begin{align}\begin{aligned}\label{e.230}
\Pi_{V_+}(f(\cdot,-\onehalf))=0,\quad \Pi_{V_-}(f(\cdot,\onehalf))=0,
\end{aligned}\end{align}
where the function $F_T$ is defined in \eqref{e.27}. As in Figure \ref{figure9}, we set
\begin{align}\begin{aligned}\label{e.311}
\widetilde{N}=Y_{[-\onehalf,\frac{2}{T}]}\cong N\backslash Y_{[\frac{2}{T},\onehalf]}.
\end{aligned}\end{align}
By the same arguments as in Proposition \ref{p.1}, for $m>0$ sufficiently large we have
\begin{align}\begin{aligned}\label{e.310}
&\eta((D_{N}+mF_T\GS)_{V_+,V_-})-\eta((D_{N}+m\GS)_{V_+,V_-})\\
=&\eta((D_{\widetilde{N}}+mF_T\GS)_{V_+,V_-})-\eta((D_{\widetilde{N}}+m\GS)_{V_+,V_-}),
\end{aligned}\end{align}
since $D_N+mF_T\GS$, $D_N+m\GS$ coincide when restricted to $Y_{[\frac{2}{T},\onehalf]}$.
For $T>0$ sufficiently large, the operator $(D_{\widetilde{N}}+mF_T\GS)_{V_+,V_-}$ is a small perturbation
of $(D_{\widetilde{N}}-m\GS)_{V_+,V_-}$ on $\widetilde{N}$, hence the associated eta invariants coincide, i.e.,
\begin{align}\begin{aligned}\label{e.312}
\eta((D_{\widetilde{N}}+mF_T\GS)_{V_+,V_-})=\eta((D_{\widetilde{N}}-m\GS)_{V_+,V_-}).
\end{aligned}\end{align}

\begin{thm}\label{t.9}
  For $m>0$, $T>0$ sufficiently large, we have
\begin{align}\begin{aligned}\label{e.240}
&\eta((D_N+mF_T\GS)_{V_+,V_-})-\eta((D_{N}+m\GS)_{V_+,V_-})
=-2\ind (D_N)_{V_+,V_-}.
\end{aligned}\end{align}
\end{thm}
\begin{proof}
 By \eqref{e.310} and \eqref{e.312}, we get
\begin{align}\begin{aligned}\label{e.313}
&\eta((D_{N}+mF_T\GS)_{V_+,V_-})-\eta((D_{N}+m\GS)_{V_+,V_-})\\
=&\eta((D_{\widetilde{N}}-m\GS)_{V_+,V_-})-\eta((D_{\widetilde{N}}+m\GS)_{V_+,V_-}).
\end{aligned}\end{align}
 For $m>0$, by Theorem \ref{t.2} we have
\begin{align}\begin{aligned}\label{e.314}
\eta((D_{\widetilde{N}}+m\GS)_{V_+,V_-})-\eta((D_{\widetilde{N}}-m\GS)_{V_+,V_-})=2\ind (D_{\widetilde{N}})_{V_+,V_-}.
\end{aligned}\end{align}
Since $N=\widetilde{N}\cup_Y Y_{[\frac{2}{T},\onehalf]}$ and the Fredholm index is homotopy stable when deforming the Riemannian metric, we have
\begin{align}\begin{aligned}\label{e.315}
\ind (D_{\widetilde{N}})_{V_+,V_-}=\ind (D_{N})_{V_+,V_-}.
\end{aligned}\end{align}
Equation \eqref{e.240} follows from \eqref{e.313}, \eqref{e.314} and \eqref{e.315}. The proof of Theorem \ref{t.9} is completed.
\end{proof}

\subsection{The index of $(D_N)_{V_+,V_-}$}
We try to compute the index of the Dirac operator $(D_N)_{V_+,V_-}$ over the finite cylinder $N=Y\times [-\onehalf,\onehalf]$.

Let $\phi_{\lambda_k}$ be the eigensection of $D_Y$ associated to the eigenvalue $\lambda_k>0$, i,e,
\begin{align}\begin{aligned}\label{e.250}
D_Y\phi_{\lambda_k}=\lambda_k\phi_{\lambda_k},\quad \norm{\phi_{\lambda_k}}_Y=1.
\end{aligned}\end{align}
By \eqref{e.39} and \eqref{e.250}, we have
\begin{align}\begin{aligned}\label{e.251}
D_Y\gamma\phi_{\lambda_k}=-\lambda_k\gamma\phi_{\lambda_k},\quad \norm{\gamma\phi_{\lambda_k}}_Y=1.
\end{aligned}\end{align}
Let $\set{\psi_i|1\leq i\leq n_+}$ be an orthonormal basis of $V_+$, then
\begin{align}\begin{aligned}\label{e.252}
\set{\psi_i,\,\gamma\psi_i,\,\phi_{\lambda_k},\,\gamma\phi_{\lambda_k}|\quad 1\leq i\leq n_+,\, k\in \N}
\end{aligned}\end{align}
form an orthonormal basis of $L^2(Y,S|_Y)$.
We expand the section $f\in \Cinf{N,S}$ in terms of the basis \eqref{e.252} as follows:
\begin{align}\begin{aligned}\label{e.218}
f(y,u)=\sum_{k=1}^\infty \left(g_{\lambda_k}(u)\phi_{\lambda_k}+f_{\lambda_k}(u)\gamma\phi_{\lambda_k}\right)
+\sum_{i=1}^{n_+}\left( h_i(u)\psi_i+ e_i(u)\gamma\psi_i\right).
\end{aligned}\end{align}
We impose the following boundary conditions for the Dirac operator $D_{N}$ over $N$,
\begin{align}\begin{aligned}\label{e.219}
\Pi_{V_+}(f(\cdot,-\onehalf))=0,\quad \Pi_{V_-}(f(\cdot,\onehalf))=0.
\end{aligned}\end{align}

Let $E_k={\rm span}_{\C}\left(\phi_{\lambda_k},\gamma\phi_{\lambda_k}\right)$ be the complex vector space spanned by the orthonormal eigensections $\phi_{\lambda_k},\gamma\phi_{\lambda_k}$. Set $V_i={\rm span}_{\C}\left(\psi_{i},\gamma\psi_{i}\right)$ for $1\leq i\leq n_+$. The solution of $D_Ns=\mu s$ for $s\in L^2([-\onehalf,\onehalf],E_k)$ has the following form,
\begin{align}\begin{aligned}\label{e.221}
s(u,y)=&c_1e^{u\sqrt{\lambda_k^2-\mu^2}}\left((\lambda_k-\sqrt{\lambda_k^2-\mu^2})\phi_{\lambda_k}(y)+\mu\gamma\phi_{\lambda_k}(y)\right)\\
&+c_2e^{-u\sqrt{\lambda_k^2-\mu^2}}\left(\mu\phi_{\lambda_k}(y)+(\lambda_k-\sqrt{\lambda_k^2-\mu^2})\gamma\phi_{\lambda_k}(y)\right).
\end{aligned}\end{align}
By the boundary condition \eqref{e.219},  we have
  \begin{align}\begin{aligned}\label{e.223}
\mu^2=\lambda_k^2+\pi^2j^2,\quad \text{for}\quad j\in \Z,\,\lambda_k>0.
\end{aligned}\end{align}

The solution of $D_Ns=\mu s$ for $s\in L^2([-\onehalf,\onehalf],V_i)$ has the following form,
\begin{align}\begin{aligned}\label{e.253}
s(u,y)=a_1\left(\cos (\mu u)\psi_i(y)-\sin(\mu u)\gamma\psi_i(y)\right)+a_2\left(\sin(\mu u)\psi_i(y)+\cos(\mu u)\gamma\psi_i(y)\right).
\end{aligned}\end{align}
By the boundary condition \eqref{e.219},  we get
\begin{align}\begin{aligned}\label{e.254}
\mu=\pi j,\quad \text{for}\quad j\in \Z.
\end{aligned}\end{align}
By \eqref{e.223}, \eqref{e.253} and \eqref{e.254}, the kernel of $(D_N)_{V_+,V_-}$ is spanned by
\begin{align}\begin{aligned}\label{e.255}
\left\{\gamma\psi_1,\gamma\psi_2,\cdots,\gamma\psi_{n_+}\right\}.
\end{aligned}\end{align}
By \eqref{e.255}, we have
\begin{align}\begin{aligned}\label{e.279}
\ind (D_N)_{V_+,V_-}&=\tr(\GS|_{\Ker (D_N)_{V_+,V_-}})\\
&=\tr(\GS|_{V_-})=-n_+,
\end{aligned}\end{align}
where $n_+=\dim V_+$.

\subsection{Proof of Theorem \ref{t.3}}

First we try to handle the following term appearing in \eqref{e.217},
\begin{align}\begin{aligned}\label{e.244}
\eta((D_{M_2}+m\GS)_{V_+})-\eta((D_{M_2}-m\GS)_{V_+}).
\end{aligned}\end{align}

For $m>0$, by Theorem \ref{t.2} we have
\begin{align}\begin{aligned}\label{e.277}
\eta((D_{M_2}+m\GS)_{V_+})-\eta((D_{M_2}-m\GS)_{V_+})=2\ind (D_{M_2,V_+}).
\end{aligned}\end{align}

Since $M_2=Y_{[-\onehalf,0]}\cup _Y M$ and the index of Fredholm operator is homotopy stable, we have
\begin{align}\begin{aligned}\label{e.245}
\ind (D_{M_2,V_+})=\ind (D_{M,V_+}).
\end{aligned}\end{align}
For $m>0$ sufficiently large and the Lagrangian subspace $V_+\subset \Ker D_Y$, by \eqref{e.277} and \eqref{e.245} we have
\begin{align}\begin{aligned}\label{e.246}
\ind D_{M,V_+}=\frac{\eta((D_{M_2}+m\GS)_{V_+})-\eta((D_{M_2}-m\GS)_{V_+})}{2}.
\end{aligned}\end{align}

\begin{proof}[Proof of Theorem \ref{t.3}]
By Theorem \ref{t.9}, \eqref{e.32}, \eqref{e.217}, \eqref{e.279} and \eqref{e.246}, for $T>0$ and $m>0$ sufficiently large we get
\begin{align}\begin{aligned}\label{e.278}
&\eta(D_{\bM}+m\kappa\GS)-\eta(D_{\bM}-m\GS)\\
=&\eta(D_{\bM}+mF_T\GS)-\eta(D_{\bM}-m\GS)\\
=&\eta((D_{M_2}+m\GS)_{V_+})-\eta((D_{M_2}-m\GS)_{V_+})\\
&\quad\quad +\eta((D_{N}+mF_T\GS)_{V_+,V_-})-\eta((D_{N}+m\GS)_{V_+,V_-})\\
=&2\ind D_{M,V_+}-2\ind (D_N)_{V_+,V_-}\\
=&2\ind D_{M,V_+}+2n_+,
\end{aligned}\end{align}
where $n_+=\dim V_+$.
By \eqref{e.278}, we get
\begin{align}\begin{aligned}\label{e.249}
\ind D_{M,V_+}=\frac{\eta(D_{\bM}+m\kappa\GS)-\eta(D_{\bM}-m\GS)}{2}-n_+.
\end{aligned}\end{align}
The proof of Theorem \ref{t.3} is completed.
\end{proof}

\subsection{Proof of Theorem \ref{t.5}}
To prove Theorem \ref{t.5}, we only need to compare the indices of the Dirac operators $D_{M,V_+}$ and $D_{M,P_{\geq}}$. By \cite[Def. 15.8]{BW93}, we introduce the concept of virtual codimension.

\begin{defn}
Let $P_1,\,P_2$ be pseudo-differential projections with the same principal symbol. The integer
\begin{align}\begin{aligned}\label{e.282}
\mathbf{i}(P_2,P_1):=\ind \left\{P_2P_1:\ran(P_1)\ra \ran(P_2)\right\}
\end{aligned}\end{align}
is called the virtual codimension of $P_2$ in $P_1$.
\end{defn}

The operator $P_2P_1:\ran(P_1)\ra \ran(P_2)$ is Fredholm, since
\begin{align}\begin{aligned}\label{e.284}
\Ker\set{P_2P_1:\ran P_1\ra \ran P_2}=\Ker T\cap \ran P_1,
\end{aligned}\end{align}
and
\begin{align}\begin{aligned}\label{e.285}
\ran \set{P_2P_1:\ran P_1\ra \ran P_2}=\ran T\cap \ran P_2,
\end{aligned}\end{align}
where the $0-$order pseudo-differential operator
\begin{align}\begin{aligned}\label{e.286}
T:=P_2P_1+(\Id-P_2)(\Id-P_1)
\end{aligned}\end{align}
has principal symbol equal $\Id$, hence is elliptic.

\begin{lemma}\label{l.2} For the projection operators $P_{\geq}$ and $\Pi_{V_+}$ (see \eqref{e.196}, \eqref{e.124}), we have
\begin{align}\begin{aligned}\label{e.287}
\mathbf{i}(\Pi_{V_+},P_{\geq})=n_+.
\end{aligned}\end{align}
\end{lemma}
\begin{proof}
For the projection operators $P_{\geq}$ and $\Pi_{V_+}$, by \eqref{e.196}, \eqref{e.124} and \eqref{e.286} we have
  \begin{align}\begin{aligned}\label{e.288}
T:&=\Pi_{V_+}P_{\geq}+(\Id-\Pi_{V_+})(\Id-P_{\geq})\\
&=\Pi_{V_+}+\Pi_{V_-}P_{<}=\Pi_{V_+}+P_{<}\\
&=\Id-\pr{V_-}.
\end{aligned}\end{align}
By \eqref{e.288}, we get
 \begin{align}\begin{aligned}\label{e.289}
\Ker T=V_-,\quad \ran T=\ran \Pi_{V_+}\oplus \ran P_<.
\end{aligned}\end{align}
By \eqref{e.284}, \eqref{e.285} and \eqref{e.288}, we have
\begin{align}\begin{aligned}\label{e.290}
&\Ker\set{\Pi_{V_+}P_{\geq}:\ran P_{\geq}\ra \ran \Pi_{V_+}}\\
=&\Ker T\cap \ran P_{\geq}
=V_-\cap \ran P_{\geq}=V_-,\\
&\ran \set{\Pi_{V_+}P_{\geq}:\ran P_{\geq}\ra \ran \Pi_{V_+}}\\
=&\ran T\cap \ran \Pi_{V_+}=\ran \Pi_{V_+}.
\end{aligned}\end{align}
By \eqref{e.282} and \eqref{e.290}, we get
\begin{align}\begin{aligned}\label{e.291}
\mathbf{i}(\Pi_{V_+},P_{\geq})=&\dim \Ker \set{\Pi_{V_+}P_{\geq}:\ran P_{\geq}\ra \ran \Pi_{V_+}}\\
&-\dim \Coker \set{\Pi_{V_+}P_{\geq}:\ran P_{\geq}\ra \ran \Pi_{V_+}}\\
=&\dim V_--\dim \frac{\ran \Pi_{V_+}}{\ran \Pi_{V_+}}=n_+.
\end{aligned}\end{align}
The proof of Lemma \ref{l.2} is completed.
\end{proof}

\begin{proof}[Proof of Theorem \ref{t.5}]
By Proposition 21.4 in \cite{BW93} and Lemma \ref{l.2}, we have
\begin{align}\begin{aligned}\label{e.281}
\ind D_{M,V_+}-\ind D_{M,P_{\geq}}=\mathbf{i}(\Pi_{V_+},P_{\geq})=n_+.
\end{aligned}\end{align}
By \eqref{e.249} and \eqref{e.281}, we get
\begin{align}\begin{aligned}\label{e.283}
\ind D_{M,P_{\geq}}=\frac{\eta(D_{\bM}+m\kappa\GS)-\eta(D_{\bM}-m\GS)}{2}-2n_+.
\end{aligned}\end{align}
The proof of Theorem \ref{t.5} is completed.
\end{proof}


\end{document}